\documentclass[a4paper]{amsart}
\usepackage{amsmath}
\usepackage{amsthm}
\usepackage{amsfonts,amssymb,mathrsfs}
\usepackage{mathtools}
\usepackage[a4paper, hmargin=3.5cm, vmargin=2.8cm]{geometry}
\usepackage{enumitem}
\usepackage{hyperref}
\usepackage{aligned-overset}
\usepackage[svgnames]{xcolor}

\usepackage{graphicx}
\usepackage{tikz}
\usetikzlibrary{calc,arrows.meta,decorations.pathreplacing,intersections}

\title[$L^1$ Localization for extended-real-valued Lagrangians]{$L^1$ Localization for extended-real-valued Lagrangians and the electromagnetic Monge problem}
\author{Rotem Assouline}
\address{Department of Mathematics, Massachusetts Institute of Technology, Cambridge, MA 02139, USA}
\email{rotem\_a@mit.edu}
\author{Davide Manini}
\address{Fakultät für Mathematik, Universität Wien, Oskar-Morgenstern-Platz 1, 1010 Wien, Austria}
\email{davide.manini@univie.ac.at}
\date{November 27, 2026}

\newenvironment{declaration}[1]{%
  \setlength{\parindent}{0em}%
  \setlength{\parskip}{-0.0mm}%
  \par%
  \vskip 1mm\footnotesize\textbf{#1.}\;\;}{\par}

\setlist[enumerate]{label={(\roman*)}}

\newcommand{\pushf}{\sharp}
\newcommand{\I}{\mathcal{I}}
\newcommand{\mmz}{\mathscr{G}}
\newcommand{\mm}{\mathfrak{m}}
\newcommand{\Leb}{\mathscr{L}}
\newcommand{\interior}{\mathrm{int}}
\newcommand{\inter}{\lambda} 

\DeclareMathOperator{\dom}{dom}
\newcommand{\eps}{\varepsilon}
\newcommand{\vphi}{\varphi}
\DeclareMathOperator{\Len}{Len}
\DeclareMathOperator{\id}{Id}
\renewcommand{\div}{{ \rm div}}
\DeclareMathOperator{\supp}{supp}

\DeclareMathOperator{\Ric}{Ric}

\newcommand{\de}{\mathrm{d}}

\newcommand{\RR}{\mathbb{R}}
\newcommand{\NN}{\mathbb{N}}

\newcommand{\CC}{\mathrm{C}}
\newcommand{\rT}{\mathrm{T}}
\newcommand{\cc}{\mathrm{c}}
\newcommand{\action}{\mathrm{S}}

\newcommand{\cL}{\mathcal{L}}
\newcommand{\cA}{\mathcal{A}}
\newcommand{\cT}{\mathcal{T}}
\newcommand{\cX}{\mathcal{X}}

\newcommand{\bz}{\mathbf{0}}

\newcommand{\sW}{\mathscr{W}}

\newcommand{\rayproj}{R^{i,k}}

\newtheorem{theorem}{Theorem}[section]
\newtheorem{lemma}[theorem]{Lemma}
\newtheorem{proposition}[theorem]{Proposition}

\theoremstyle{definition}
\newtheorem{definition}[theorem]{Definition}
\newtheorem{remark}[theorem]{Remark}
\newtheorem*{remark*}{Remark}

\begin{document}

\begin{abstract}
    We construct $L^1$ displacement interpolations for electromagnetic Lagrangians
    on globally-hyperbolic Lorentzian manifolds.
    Our construction relies on an extension of Klartag's needle decomposition technique to weighted manifolds endowed with
    a Lagrangian which is smooth and strictly convex on its domain,
    but which may attain the value $\infty$ and may not grow superlinearly on each fiber.
\end{abstract}

\maketitle

    \section{Introduction}

    The solution to the optimal transport problem with cost arising from a Tonelli Lagrangian 
    gives rise to a one-parameter family of measures interpolating
    between a pair of given probability measures~\cite{BernardBuffoni,FathiFigalli,Villani09}. 
    In the case of $L^1$ (Monge{--}Kantorovich) optimal transport,
    where the cost is the free-time action potential,
    more is true: the manifold decomposes
    into disjoint action-minimizing trajectories,
    and by transporting mass along each such trajectory
    one obtains an optimal transport map~\cite{EvGa99,CFM,FeMc02,BernardBuffoni06,Figalli07}.
    This property of the Monge-Kantorovich problem can be used
    to reduce geometric inequalities on weighted manifolds to their one-dimensional counterparts via a method
    known as \emph{needle decomposition}, due originally to Klartag~\cite{Kl17}.
    Klartag's needle decomposition technique, which was inspired by the localization technique in convex geometry~\cite{GrMi87,KaLoSi95},
    has since been extended to a variety of settings~\cite{CavallettiMondino17,Ohta18,Milman21,Assouline24}.
    
    \medskip
    Lower bounds on the Ricci curvature of Riemannian metrics can be characterized in terms
    of convexity properties of entropy functionals along displacement interpolations
    for the squared distance cost~\cite{CorderoErausquinMcCannSchmuckenschlaeger01,RenesseSturm05,St06-2,St06,LoVi09}
    as well as in terms of corresponding bounds on the individual needles in the needle decomposition~\cite{Kl17}.
    Analogous results in other metric settings were established in~\cite{Ohta09a,McCann20,CavallettiMondino24,CaMi21,BraunMcCann24,BraunOhta24,BaKrSi18,BaRi19,BarilariMondinoRizzi26,Milman21}.
    For general Tonelli Lagrangians, one can define a quantity which generalizes Ricci curvature
    and which plays a similar role in determining the behavior of entropy along displacement interpolations~\cite{lee2013displacement,Ohta+2014,Assouline24}.
    
    \medskip
    In Lorentzian geometry,
    where trajectories are constrained to be future-directed, 
    the optimal transport problem can be formulated using a Lagrangian which is infinite outside the future cone in each tangent space~\cite{McCann20,BraunGigliMcCannOhanyanSaemann26}.
    A natural choice is the Lorentzian norm raised to a power $r \in (0,1)$;
    besides attaining the value $\infty$, this Lagrangian fails to grow superlinearly at infinity,
    a condition required of Tonelli Lagrangians.
    
    \medskip
    In this work we consider Lagrangians which may exhibit the above departures from the Tonelli setting,
    and construct displacement interpolations and needle decompositions for them.
    The main application we have in mind is optimal transport with cost arising from an
    electromagnetic Lagrangian combining a Lorentzian metric with a linear term representing the electromagnetic potential.
    In a companion work~\cite{AssoulineManini26}, we use these results
    to establish a connection between displacement convexity
    along optimal transport with electromagnetic cost
    and lower bounds on a suitably defined generalization
    of the Ricci curvature.
    We hope that the technical groundwork laid here will be useful for studying 
    optimal transport with other Lagrangian costs. 
    
    \medskip
    A major difficulty in this setting is 
    to determine under which conditions displacement interpolations exist at all.
    In the electromagnetic case, the works of Minguzzi and Sánchez~\cite{Minguzzi03,Minguzzi04,MinguzziSanchez06} 
    on the existence of action-minimizing solutions of the Lorentz force equation joining given pairs of points
    are of great importance in this regard, see Section~\ref{sec:em}.

    \medskip
    We remark that Lagrangians arising in optimal control theory
    do not fall within the scope of this work,
    as their domain meets each tangent space in a proper subspace,
    whereas we require the interior of the domain to intersect every fiber.
    For optimal transport in the setting of control theory
    see~\cite{AgrachevLee09}; see~\cite{FigalliRifford10,BorzaKlingenbergWood25,BorzaRigoniZoghlami25}
    for the particular cases of sub-Riemannian and sub-Lorentzian
    geometry.

    \subsection{Displacement interpolations}

    Let $M$ be a smooth manifold and let $L : TM \to \RR\cup\{+\infty\}$ be a Lagrangian, which we assume to be convex and lower semicontinuous on each tangent space,
    smooth with positive definite vertical Hessian on the interior of
    its domain, positive on the zero section, and essentially smooth
    (see Assumptions~\ref{item:convex}{--}%
    \ref{item:positivezerosection} in Section~\ref{sec:Lagrangians}).
    We denote by $H:T^*M \to \RR\cup\{+\infty\}$ the \emph{Hamiltonian}
    associated to $L$, and by $\cL$ the \emph{Legendre transform},
    which is a diffeomorphism between the interiors of the domains of $H$ and $L$.
    
    \medskip
    For $x,y \in M$, the cost
    \[
      \cc(x,y)
      =
      \inf_\gamma\int L(\dot\gamma)
    \]
    is defined to be the infimum of the action over piecewise-$C^1$ curves $\gamma$ joining $x$ to
    $y$, defined on an arbitrary finite interval, whose velocity lies in the interior of the domain of $L$.
    A curve attaining this infimum is called a \emph{minimizing extremal}.
    We fix a set $\cA \subseteq (M\times M) \setminus \Delta$, where $\Delta$ denotes the diagonal, and assume
    that every pair of points in $\cA$ is joined by a minimizing extremal.
    We also assume that minimizing extremals whose endpoints range over a compact
    subset of $\cA$ are defined on intervals of uniformly bounded length and have
    velocities in a fixed compact subset of the interior of the domain of $L$;
    see Assumptions~\ref{item:attained} and~\ref{item:compactness}
    in Section~\ref{sec:ACM}.
    Sufficient conditions on $\cA$ that ensure those properties in the electromagnetic case
    are provided in Section~\ref{sec:em}.
  
    \medskip
    Our first result asserts the existence of a dynamical solution to the optimal transport problem in the following sense, cf.~\cite{BernardBuffoni06,FathiFigalli,Villani09,Assouline24}.
    Let $\Phi$ denote the Euler{--}Lagrange flow on the tangent bundle, so that
    for every tangent vector $v$ in the interior of the domain of $L$,
    the curve $\tau \mapsto \pi(\Phi(v,\tau))$ is the solution of the
    Euler{--}Lagrange equation with initial velocity $v$, where
    $\pi : TM \to M$ is the projection.
    A \emph{displacement interpolation} between a pair of probability measures $\mu_0,\mu_1$ on $M$ is a family $(\mu_\inter)_{\inter\in[0,1]}$ of probability measures of the form
    \[
      \mu_\inter = (\rT_\inter)_\sharp\mu_0, \qquad
      \rT_\inter := \pi\circ\Phi\circ(V_0,\inter\ell), \qquad 0 \le \inter \le 1,
    \]
    where $V_0$ is a vector field on $\supp\mu_0$ and $\ell \ge 0$ is a function on $\supp\mu_0$, such that for $\mu_0$-almost every $x_0$, the curve $\tau \mapsto \pi\left(\Phi(V_0\vert_{x_0},\tau)\right),\,\,0 \le \tau \le \ell(x_0)$, is a minimizing extremal, and the map $\rT_1$ is an optimal transport map from $\mu_0$ to $\mu_1$ with respect to the cost $\cc$, i.e., minimizes $\int_M \cc(x_0,\rT(x_0))\,\de\mu_0(x_0)$ among all maps pushing forward $\mu_0$ to $\mu_1$.
    A displacement interpolation is called \emph{monotone} if $\Phi(V_0\vert_{x_0},\tau) = V_0\vert_{x_0'}$ implies $\ell(x_0) \le \ell(x_0') + \tau$, for all $x_0,x_0' \in \supp\mu_0$ and all $\tau \ge 0$.
    The precise definitions are given in Section~\ref{sec:di}.

    \begin{theorem}[Existence and uniqueness of displacement interpolations]\label{thm:di-intro}
        Let $M$ be a smooth manifold, let $L:TM\to\RR\cup\{+\infty\}$ be a Lagrangian satisfying Assumptions~\ref{item:convex}{--}\ref{item:positivezerosection} of Section~\ref{sec:Lagrangians},
        and let $\cA \subseteq (M\times M) \setminus \Delta$
        satisfy Assumptions~\ref{item:attained} and~\ref{item:compactness} of Section~\ref{sec:ACM}.
        Let $\mu_0$ and $\mu_1$ be compactly supported, absolutely continuous Borel probability measures on $M$ such that $\supp\mu_0\times\supp\mu_1 \subseteq \cA$.
        Then there exists a unique monotone displacement interpolation between $\mu_0$ and $\mu_1$.
    \end{theorem}

    Theorem~\ref{thm:di-intro} is proved in Section~\ref{sec:di}, see Theorem~\ref{thm:diexistence}.

    \subsection{Application to the electromagnetic Lagrangian}

    Our next result is an application of Theorem~\ref{thm:di-intro} to the case of electromagnetic Lagrangians on Lorentzian manifolds.
    Let $(M,g)$ be a globally hyperbolic, time-oriented Lorentzian manifold.
    For $x,y \in M$, write $x \ll y$ if there is a smooth future-directed timelike curve joining $x$ to $y$.
    Let $r<1$, $r\neq0$, let $A$ be a one-form on $M$, and define a Lagrangian $L_A:TM\to\RR\cup\{+\infty\}$ by
        \[
            L_A(v) : =
            \begin{cases}
                -\tfrac1r\langle v,v\rangle^{r/2} + \tfrac{1-r}{r} - A(v), & \text{ if $v$ is future-directed or $v = 0$,}
                \\
                \infty, & \text{ otherwise.}
            \end{cases}
        \]
        If $r<0$ then we also set $L_A(v):=\infty$ if $v$ is null or $v=0$.

    \begin{theorem}\label{thm:em}
        Let $\mu_0$ and $\mu_1$ be compactly supported, absolutely continuous Borel probability measures on $M$ such that $x_0 \ll x_1$ for all $(x_0,x_1) \in \supp\mu_0\times\supp\mu_1$.
        Then there exists $s_0 > 0$, depending on $A$ and on the
        supports of $\mu_0$ and $\mu_1$, such that for every $s\in (-s_0, s_0)$ there exists a displacement interpolation between $\mu_0$ and $\mu_1$ for the cost associated with $L_{sA}$.
    \end{theorem}

    For Dirac measures, Theorem~\ref{thm:em} follows from
    a result of Caponio and Minguzzi~\cite{CaponioMinguzzi04},
    which also implies that, when no pair of points $(x_0,x_1) \in \supp\mu_0\times\supp\mu_1$
    is joined by a null geodesic, 
    and in particular when $(M,g)$ is Minkowski spacetime, one can take $s_0 = \infty$.
    More refined conditions, based on the results in~\cite{MinguzziSanchez06}, are described in Section~\ref{sec:em}.
    Theorem~\ref{thm:em} can be viewed as asserting that displacement interpolations between two fixed measures exist whenever the charge-to-mass ratio of the transported particles is sufficiently small.
    Indeed, multiplying the potential $A$ by a scalar $s$ has the same effect on the
    equations of motion as multiplying the charge-to-mass ratio of the test particles by $s$.

    \subsection{Needle decomposition}

    In order to prove Theorem~\ref{thm:di-intro},
    we establish Theorem~\ref{thm:main},
    which extends Klartag's needle decomposition technique to the setting of Theorem~\ref{thm:di-intro}.
    Theorem~\ref{thm:main} provides a disintegration of the reference measure $\mm$
    into a family of measures supported on minimizing extremals
    and satisfying a mass-balance condition.
    Clause~\ref{item:HJneedle} of Theorem~\ref{thm:main},
    which gives a description of the densities of the needles
    in terms of solutions to the Hamilton-Jacobi equation,
    can be used in tandem with Bochner's inequality
    in order to deduce convexity of entropy along
    displacement interpolations from bounds on a Ricci-type
    quantity;
    in the electromagnetic setting, this is carried out in the companion work~\cite{AssoulineManini26}.

    \medskip    
    The idea of the proof of Theorem~\ref{thm:main} is to employ a regularization procedure in order to reduce it to the needle decomposition result for Tonelli Lagrangians established in~\cite{Assouline24}.
    We first replace $L$ locally, in a neighborhood of each point of the transport set,
    by a Tonelli Lagrangian which coincides with $L$ in the direction of the transport rays.
    We then apply the needle decomposition result for Tonelli Lagrangians proved in~\cite{Assouline24}
    to each local Lagrangian,
    and glue together the local needle decompositions to obtain a global one.

    \subsection{Structure of the paper}
    In Section~\ref{sec:OT} we recall facts about optimal transport with Lagrangian cost and state our assumptions on the Lagrangian.
    In Section~\ref{sec:tonellireg},
    we describe a procedure for locally replacing the Lagrangian $L$ by a Tonelli Lagrangian.
    In Section~\ref{sec:needle},
    we combine the Tonelli regularization scheme with the needle decomposition result for Tonelli Lagrangians proved in~\cite{Assouline24} to deduce Theorem~\ref{thm:main}.
    In Section~\ref{sec:di} we prove Theorem~\ref{thm:di-intro}.
    In Section~\ref{sec:em} we specialize to the electromagnetic case and prove Theorem~\ref{thm:em}.

    \begin{declaration}{Acknowledgements}
    Part of this work was carried out while R.\ A.\ was supported by the Fondation Sciences Mathématiques de Paris (FSMP) and the Rothschild Fellowship (Yad Hanadiv).

    The research leading to these results is part of a project that
    has received funding from the European Research Council (ERC)
    under the European Union's Horizon 2020 research and innovation
    programme (grant agreement No 101001677).

    D.\ M.\ is supported by the Austrian Science Fund (FWF) [grant
    DOI \href{https://doi.org/10.55776/EFP6}{10.55776/EFP6}].
    \end{declaration}
    \begin{declaration}{Use of AI tools}
    The authors started working on this project in summer 2025.
    From spring 2026 onward the authors used AI agents, first for
    editorial tasks (e.g., bulk changes of notation) and later to check proofs and to suggest simplifications.
    Most of the mathematical content was in place before the dramatic improvement in mathematical capabilities of AI tools during summer-fall 2026.
    \end{declaration}
    \begin{declaration}{License}
    For open-access purposes, this work is licensed under a
    \href{https://creativecommons.org/licenses/by/4.0/deed.en}{``CC BY
      4.0''} license.
    \end{declaration}

    \section{Optimal transport with Lagrangian cost}\label{sec:OT}

    \subsection{Extended-real-valued Lagrangians}\label{sec:Lagrangians}

        Let $M$ be a smooth manifold and let 
        \[
          L:TM \to \RR\cup\{+\infty\}
        \]
        be a function, hereinafter referred to as the \emph{Lagrangian}, satisfying the following assumptions:
        \begin{enumerate}[series=assumptions, label={(\Roman*)}]
            \item\label{item:convex}
                For all $x \in M$, the restriction of $L$ to $T_xM$ is
                a proper (i.e., not identically $+\infty$) lower
                semicontinuous convex function.
            \item\label{item:smooth}
                On the interior of its domain $\dom L : = L^{-1}(\RR)$, the Lagrangian $L$ is smooth and its vertical Hessian $D^2L$ is positive definite.
            \item\label{item:essentiallysmooth}
                The Lagrangian $L$ is fiberwise essentially smooth in the sense of~\cite[Sect.~26]{Ro70}.
                This means that, for every $x \in M$, the set $\interior\dom L\cap
                T_xM$ is nonempty, and for every sequence $v_n \in \interior\dom L \cap T_xM$ with $v_n\to v \in \partial\dom L \cap T_xM$,
                the sequence $DL\vert_{v_n}$ tends to infinity (i.e., eventually leaves every compact set).
            \item\label{item:domaincompat}
            For every
            $x\in M$,
            \[
              \interior\dom L\cap T_xM
              =
              \interior(\dom L\cap T_xM).
            \]
            \item\label{item:positivezerosection} The Lagrangian $L$ is positive on the zero section: $0 < L\vert_\bz$.
            \end{enumerate}

        The Lagrangian $L$ is called a \emph{Tonelli Lagrangian} if it is finite everywhere, satisfies condition~\ref{item:smooth}, and is fiberwise superlinear, i.e., for every Riemannian metric $g$ on $M$ and every compact set $A \subseteq M$ there exists a constant $C > 0$ such that
        \[
          L(v) \ge |v|_g - C,
          \qquad
          \text{ for all }
          v \in T_xM, \, x \in A,
        \]
        see~\cite[Appendix~B]{FathiFigalli},~\cite[Ch.~1]{Sorrentino15} or~\cite[Sect.~2.1]{Assouline24}.

        \begin{remark*}
          Finsler metrics are not smooth at the zero section, but since we will be restricting our attention to the zeroth energy level (which we normalize to be the unit tangent bundle in the Finsler case), they can be modified to be smooth there without affecting the part of dynamics which is relevant for the optimal transport problem.
        \end{remark*}

        Under assumptions~\ref{item:convex}--\ref{item:domaincompat}, the Hamiltonian
        \begin{equation}\label{eq:Hdef}
            H: T^*M \to \RR\cup\{+\infty\}, \qquad H(p) : = \sup_{v\in T_xM} \left[p(v) - L(v)\right], \qquad p \in T^*_xM, x \in M
        \end{equation}
        is well-defined and fiberwise strictly convex on the interior
        of its domain $\dom H : = H^{-1}(\RR)$.
        The Legendre
        transform $\cL:\interior\dom H \to \interior\dom L$,
        sends a covector $p \in T^*_xM$ to the unique vector $v \in T_xM$ at which the supremum in~\eqref{eq:Hdef} is attained;
        it satisfies $\cL^{-1}(v) = DL\vert_v$ for $v \in \interior\dom L$.
        Since the map $\cL$ is a fiberwise diffeomorphism~\cite[Th.~26.5]{Ro70}, it is in fact a diffeomorphism between the
        interior of the domains by Assumptions~\ref{item:smooth} and~\ref{item:domaincompat}.
        We also have that $D^2H = (D^2L)^{-1}\circ\cL$ on
        $\interior\dom H$.
        We refer the reader to~\cite{Ro70} for more details on convex
        functions and duality.
        
        \subsection{Action, cost, and minimizing extremals}\label{sec:ACM}

        \leavevmode
        The \emph{action} of a piecewise-$C^1$ curve $\gamma:[0,\ell]\to M$ is defined by
        \[
            \action[\gamma] : = \int_0^\ell L(\dot\gamma(t))\,\de t.
        \]

        We define a cost function $\cc$ on $M\times M$ by
        \begin{equation}\label{eq:cdef}
          \cc(x,y) : = \inf_\gamma \action[\gamma]
          \in [-\infty, \infty],
        \end{equation}
        where the infimum is taken over all piecewise-$C^1$ curves
        $\gamma$ joining $x$ to $y$ and satisfying
        \[
          \dot\gamma(t)\in \interior\dom L \qquad \text{for all } t,
        \]
        with the
        convention $\inf\varnothing=+\infty$.
        Since the domain of $\gamma$ is a free
        interval, the cost $\cc$ satisfies the triangle
        inequality:
        \begin{equation}\label{eq:triangleineq}
          \cc(x,x'') \le \cc(x,x') + \cc(x',x''),
          \qquad x,x',x'' \in M,
        \end{equation}
        with the convention $+\infty-\infty = -\infty + \infty = +\infty$.
        A curve $\gamma$ attaining the infimum in~\eqref{eq:cdef} will be called a \emph{minimizing extremal}. A curve defined on an arbitrary interval is called a minimizing extremal if its restriction to every compact subinterval is.
            
        \medskip
        The \emph{energy} is the function $E:\interior\dom L \to \RR$ defined by
        \[
            E : = H \circ \cL^{-1},
        \]
        and we write
        \[
              SM : = E^{-1}(0) = \{v \in \interior\dom L \,\, \colon \, \, E(v) = 0\}.
        \]

        It is well-known that minimizing extremals are smooth and satisfy the Euler{--}Lagrange equation
        \begin{equation}\label{eq:EL}
          \tfrac{d}{dt}(\partial_vL\circ\dot\gamma) = \partial_xL\circ\dot\gamma
        \end{equation}
        in any canonical local coordinate chart $(x,v)$ on $TM$,
        and that the energy is constant along minimizing extremals.
        In fact, since the infimum in the definition of $\cc$ is taken
        over curves defined on arbitrary time intervals, the energy vanishes along minimizing extremals:
        \[
          E(\dot\gamma) \equiv 0,
        \]
        as can be seen by considering variations of the form
        $\gamma_\lambda(t) : = \gamma(\lambda t)$ for $\lambda >
        0$ and using conservation of energy.

        \begin{lemma}\label{lem:zeroreg}
          The value $0$ is a regular value of the energy $E$.
        \end{lemma}
        \begin{proof}
          By~\ref{item:positivezerosection}, it suffices to show that $E$ has no critical points on $\interior\dom L\setminus\bz$. The energy is given in local coordinates by
          \[
            E(v) = DL\vert_v\cdot v - L(v),
            \qquad
            v \in \interior\dom(L).
          \]
          
          Differentiating in the radial vertical direction and
          multiplying by $v$ gives
          \[
            DE\vert_v\cdot v = D^2L\vert_v(v,v) > 0,
          \]
          whence $v$ is not a critical point of $E$.
        \end{proof}

        \begin{remark*}
          Energy levels other than $0$ can be considered by adding a
          constant to the Lagrangian;
          if Assumptions~\ref{item:convex}{--}%
          \ref{item:positivezerosection} are satisfied by the
          Lagrangian $L$ then they are satisfied by $L + e_0$ for all
          $e_0 > -\inf_{\bz}L$.
          However, the global assumptions below may not persist.
        \end{remark*}

        Assumptions~\ref{item:convex}-\ref{item:positivezerosection} are local;
        we now impose some global conditions.
        Fix a set 
        \[
          \cA \subseteq (M\times M)\setminus\Delta,
        \]
        where $\Delta \subseteq M\times M$ denotes the diagonal, with the following properties:
        \begin{enumerate}[resume*=assumptions]
            \item\label{item:attained} For every $(x_0,x_1) \in
              \cA$, there exists a minimizing extremal
              joining $x_0$ to $x_1$.
            \item\label{item:compactness} For every compact set $K \subseteq \cA$ there exist a compact set $V \subseteq \interior\dom L$ and $T > 0$ such that every minimizing extremal $\gamma:[0,\ell]\to M$ with $(\gamma(0),\gamma(\ell))\in K$ satisfies $\ell \le T$ and $\dot\gamma(t) \in V$ for all $t \in [0,\ell]$.
        \end{enumerate}

        The next lemma asserts that the cost is upper
        semicontinuous, and continuous on the set $\cA$;
        the proof of upper semicontinuity is similar to the proof of lower semicontinuity
        of the time-separation function in Lorentzian
        geometry (see, e.g,~\cite[Prop.~14.17]{ONeill83}).

        \begin{lemma}\label{lem:cost-usc}
          Under Assumptions~\ref{item:convex}{--}\ref{item:positivezerosection}, the cost $\cc$ is upper semicontinuous on $M\times M$.
          If, in addition, Assumptions~\ref{item:attained} and~\ref{item:compactness} hold, then the cost $\cc$ is real-valued and continuous on the set $\cA$.
        \end{lemma}
        \begin{proof}
          We first prove upper semicontinuity.
          Fix $(x,y)$,
          and consider three cases: 1) $\cc(x,y)=\infty$; 2)
          $\cc(x,y)\in\RR$; 3) $\cc(x,y)=-\infty$.
          In the first case there is nothing to prove.
          We consider only the second case; the third is proved similarly.

          \medskip
          Fix $\eps > 0$, and let $\gamma:[0,\ell]\to M$ be a
          piecewise $C^1$ curve joining $x$ to $y$ and satisfying
          \begin{equation}
            \label{eq:almost-extremal}
            \cc(x,y)\geq \action[\gamma]-\eps
            .
          \end{equation}
        
        Since the action of $\gamma$ is finite, there exists $\delta > 0$ such that
        \begin{equation}\label{eq:endarcs}
          \big|\action[\gamma\vert_{[0,\delta]}]\big| \le \eps
          \qquad\text{and}\qquad
          \big|\action[\gamma\vert_{[\ell-\delta,\ell]}]\big| \le \eps
          .
        \end{equation}
        
        We may assume (up to taking $\delta>0$ smaller) that
        $\gamma|_{[0,\delta]}$ and $\gamma|_{[\ell-\delta,\ell]}$ are
        $C^1$.
        Since $\action$ is continuous in the $C^1$ topology on
        curves with velocity in $\interior\dom L$, there is a
        $C^1$-neighborhood $\Omega^-$ of $\gamma\vert_{[0,\delta]}$ on
        which $\action \le 2\eps$.
        Therefore, the set
        $U^{-}$ 
        of points
        joined to $\gamma(\delta)$ by a curve of $\Omega^-$ is a
        neighborhood of $x$.
        Define $U^+ \ni y$ analogously.
        Then
        \begin{equation}
          \label{eq:neighborhoods-semicontinuity-cost}
          \cc(x',\gamma(\delta)) \le 2\eps,
          \quad
          \forall x'\in U^-
          \qquad
          \text{ and }
          \qquad
          \cc(\gamma(\ell-\delta),y') \le 2\eps,
          \quad
          \forall y'\in U^+,
        \end{equation}
        and thus for $(x',y')\in U^-\times U^+$,
        \begin{align*}
          \cc(x',y')
          \overset{\eqref{eq:triangleineq}}&{\le}
          \cc(x',\gamma(\delta))
          +
          \cc(\gamma(\delta),\gamma(\ell-\delta))
          +
          \cc(\gamma(\ell-\delta),y')
          \\
          \overset{\text{\eqref{eq:neighborhoods-semicontinuity-cost}}}
          &\leq
          2\eps
          +
          \action\big[\gamma|_{[\delta,\ell-\delta]}\big]
          +
          2\eps
            \overset{\text{\eqref{eq:endarcs}}}\le
          6\eps
          +
            \action[\gamma]
            \overset{\text{\eqref{eq:almost-extremal}}}\leq
          7\eps
          +
            \cc(x,y)
            .
        \end{align*}
        Thus $\limsup_{(x',y')\to(x,y)}\cc(x',y') \le \cc(x,y)$, and upper semicontinuity is proved.

        \medskip
          Now suppose that Assumptions~\ref{item:attained} and~\ref{item:compactness} hold.
          Let 
          \[
            \cA \ni (x_n,y_n)  \to (x,y) \in \cA.
          \]
          
          By~\ref{item:attained}, there exist minimizing extremals $\gamma_n:[0,\ell_n]\to M$ joining $x_n$ to $y_n$, and by~\ref{item:compactness} applied to $K = \{(x_n,y_n)\}_n\cup\{(x,y)\}$, there exist a compact set $V\subseteq \interior\dom L$ and $T > 0$ such that
          \[
            \ell_n \le T \qquad \text{and} \qquad \dot\gamma_n(t) \in V \qquad \text{ for all $t \in [0,\ell_n]$.}
          \]
          
          Passing to a subsequence, we may assume that $\ell_n \to \ell_\infty \in (0,T]$; the case $\ell_\infty = 0$ is ruled out since $(x,y) \in \cA$, and the set $\cA$ is assumed to be disjoint from the diagonal.
          
          \medskip
          Since $\gamma_n$ is a sequence of solutions to the Euler{--}Lagrange equation with velocities contained in a compact subset of $\interior\dom L$, the second derivatives of $\gamma_n$ are contained in a compact subset of $TTM$; hence, by the Arzel\`a{--}Ascoli theorem, after passing to a further subsequence, the curves $\eta_n(s):=\gamma_n(\ell_ns)$ converge in $C^1([0,1],M)$ to a limit curve $\eta$.
          The curve $\gamma(t):=\eta(t/\ell_\infty)$ joins $x$ to $y$ and satisfies $\action[\gamma]=\lim_n\action[\gamma_n]$, whence $\cc(x,y) \le \action[\gamma] = \lim_n\action[\gamma_n]=\lim_n\cc(x_n,y_n)$.
          Since $(x_n,y_n)$ is an arbitrary sequence in $\cA$ converging to $(x,y)$,
          it follows that $\cc$ is lower semicontinuous at $(x,y)$. 
        \end{proof}

        \begin{remark}
          Let us compare Assumptions~\ref{item:convex}{--}%
          \ref{item:compactness} to the assumptions
          in~\cite{Assouline24}, where $L$ is assumed to be Tonelli  (and in particular real-valued).
          Assumption~\ref{item:domaincompat} is trivially true for
          real-valued Lagrangians.
          Assumption~(I) of~\cite{Assouline24} consists of Assumptions~\ref{item:convex}{--}\ref{item:essentiallysmooth} and~\ref{item:positivezerosection}, together with superlinearity and supercriticality, which we do not assume.
          If $\cA= (M\times M)\setminus\Delta$,
          then Assumptions~\ref{item:attained}
          and~\ref{item:compactness} imply Assumptions~(II)
          and~(III) in~\cite{Assouline24}, respectively.
          Therefore, for Tonelli Lagrangians, our
          paper recovers the setting of~\cite{Assouline24},
          except the requirement that the measures have disjoint supports (see Remark~\ref{rem:overlapping-supports}).
        \end{remark}

        \subsection{The Monge{--}Kantorovich problem}

        In this subsection we recall the optimal transport problem
        and its dual problem.
        We refer the reader to~\cite{Villani09} for more
        details on optimal transport.
        From this point onwards, we fix a smooth manifold $M$,
        a Lagrangian $L$ satisfying Assumptions~\ref{item:convex}{--}\ref{item:positivezerosection} and a set $\cA \subseteq (M\times M)\setminus\Delta$
        satisfying Assumptions~\ref{item:attained} and~\ref{item:compactness}.

        \medskip
        Let $\mu_0,\mu_1$ be absolutely continuous probability measures on $M$
        with compact supports $\cX_0$ and $\cX_1$, respectively.
        Assume that
        \[
          \cX_0\times\cX_1 \subseteq \cA;
        \]
        in particular, the sets $\cX_0$ and $\cX_1$ are disjoint;
        this assumption is discussed in Remark~\ref{rem:overlapping-supports}.
        Set
        \begin{equation}\label{eq:CC}
            \CC(\mu_0,\mu_1) : =\inf_\Pi\int_{\cX_0\times \cX_1}\cc(x_0,x_1)\,\de\Pi(x_0,x_1),
        \end{equation}
        where the infimum is over all couplings between $\mu_0$ and $\mu_1$.
        A coupling $\Pi$ attaining this infimum is called an \emph{optimal coupling} between $\mu_0$ and $\mu_1$.
        A measurable map $\rT : M \to M$ is called an \emph{optimal transport map}
        if the measure $(\id,\rT)_\sharp\mu_0$ is an optimal coupling between
        $\mu_0$ and $\mu_1$, and in particular $\rT_\sharp\mu_0 = \mu_1$.
        
        \medskip
        The dual problem to~\eqref{eq:CC}
        (see~\cite[Ch.~5]{Villani09}) is
        \begin{equation}\label{eq:Kdef}
            \mathrm{K}(\mu_0,\mu_1) : = \sup_{(\phi_0,\phi_1)}\int_M \phi_1\,\de\mu_1 - \int_M \phi_0\,\de\mu_0,
        \end{equation}
        where the supremum is over all pairs $(\phi_0,\phi_1) \in L^1(\mu_0)\times L^1(\mu_1)$ satisfying
        \begin{equation}\label{eq:phi0phi1Kant0}
            \phi_1(x_1) - \phi_0(x_0) \le \cc(x_0,x_1), \qquad \text{for all $(x_0,x_1) \in \cX_0\times \cX_1$}.
        \end{equation}

        A pair $(\phi_0,\phi_1)$ attaining the supremum in~\eqref{eq:Kdef} is called a \emph{Kantorovich pair} for $(\mu_0,\mu_1)$.

        \begin{proposition}[Kantorovich duality]\label{prop:KD}
            We have
            \[
                \mathrm{K}(\mu_0,\mu_1) = \CC(\mu_0,\mu_1),
            \]
            and both optimization problems admit solutions.
            Furthermore, for every pair $(\phi_0,\phi_1)$ satisfying~\eqref{eq:phi0phi1Kant0} and every coupling $\Pi$ between $\mu_0$ and $\mu_1$ the following are equivalent:
            \begin{itemize}
                \item For $\Pi$-a.e. $(x_0,x_1) \in \cX_0\times \cX_1$, \[\phi_1(x_1) - \phi_0(x_0) = \cc(x_0,x_1).\]
                \item The coupling $\Pi$ is an optimal coupling between $\mu_0$ and $\mu_1$ and the pair $(\phi_0,\phi_1)$ is a Kantorovich pair for $(\mu_0,\mu_1)$.
            \end{itemize}
        \end{proposition}
        \begin{proof}
            It suffices to prove that the cost $\cc$ on
            $\cX_0\times \cX_1$ satisfies the conditions of~\cite[Th.~5.10]{Villani09}.
            By Lemma~\ref{lem:cost-usc} and the fact that $\cX_0\times\cX_1$ is compact and contained in $\cA$,
            the function $\cc$ is finite and continuous on
            $\cX_0\times \cX_1$.
            The cost $\CC(\mu_0,\mu_1)$ is finite since
            $\int \cc\,\de(\mu_0\otimes \mu_1) < \infty$.
            We also have two-sided bounds for $\cc$ on $\cX_0\times \cX_1$
            by integrable continuous functions (indeed, by constants).
        \end{proof}

        \begin{remark}\label{rem:overlapping-supports}
          The assumption that the sets $\cX_0$ and $\cX_1$ are disjoint is used only in the proof of Proposition~\ref{prop:KD} to guarantee that the cost $\cc$ is continuous on $\cX_0\times\cX_1$; elsewhere we use only mutual singularity $\mu_0\perp\mu_1$.
          Hence Theorems~\ref{thm:main} and~\ref{thm:diexistence} remain valid when $\mu_0 \perp \mu_1$, provided that the cost $\cc$ is continuous on $\cX_0\times\cX_1$ and that Assumptions~\ref{item:attained} and~\ref{item:compactness} are
          satisfied by $\cA\setminus\Delta$ and~$\cA$, respectively (for instance, in the Riemannian setting).
          If, in addition, the cost $\cc$ is continuous on $(\cX_0\cup\cX_1)\times(\cX_0\cup\cX_1)$ and vanishes on the diagonal, then, in Theorem~\ref{thm:main}, the assumption of mutual singularity can be removed as well.
          Indeed, the optimal cost equals the maximum of $\int u\,\de(\mu_1-\mu_0)$ over all $L$-dominated functions $u$, see~Lemma~\ref{lem:strain}\ref{item:u-facts} and~\ref{item:ends} (see also~\cite[Particular Case~5.4]{Villani09}),
          and an optimal coupling between the positive and negative parts $(\mu_0-\mu_1)^\pm$ of the signed measure $\mu_0-\mu_1$
          becomes an optimal coupling between $\mu_0$ and $\mu_1$ upon adding $(\id,\id)_\sharp\big(\mu_0-(\mu_0-\mu_1)^+\big)$.
        \end{remark}

        Let $(\phi_0,\phi_1) \in L^1(\mu_0)\times L^1(\mu_1)$ be a Kantorovich pair for $\mu_0$ and $\mu_1$.
        We may assume without loss of generality (see~\cite[Th.~5.10]{Villani09}) that
        \begin{equation}\label{eq:phi0phi1normalization}
          \phi_1 =
          \min_{x_0 \in \cX_0}\left[\phi_0(x_0) +
            \cc(x_0,\cdot)\right]
          \qquad \text{and} \qquad
          \phi_0 =
          \max_{x_1\in\cX_1}\left[\phi_1(x_1) -
            \cc(\cdot,x_1)\right].
        \end{equation}
        
        In particular, $\phi_0$ and $\phi_1$ are both bounded and
        uniformly continuous (as the infimum, resp.\ supremum,
        of a family of functions with a common modulus of continuity).

        \begin{lemma}\label{lem:calibrated}
          Let $\eta:[t_0,t_1]\to M$ be a piecewise-$C^1$ curve such that
          \begin{equation}\label{eq:Kantequality}
            \eta(t_0)\in\cX_0,
            \quad
            \eta(t_1)\in\cX_1,
            \qquad
            \phi_1(\eta(t_1)) - \phi_0(\eta(t_0))
            \geq
            \action[\eta]
          \end{equation}
          and $\dot\eta(t)\in\interior\dom L$ for all $t\in[t_0,t_1]$.
          Then $\eta$ is a minimizing extremal, and there exists a maximally defined curve $\bar\gamma:[0,\ell]\to M$ extending $\eta$ and satisfying~\eqref{eq:Kantequality}.
        \end{lemma}
        \begin{proof}
          Since
          \begin{equation}
            \cc(\eta(t_0),\eta(t_1))
            \overset{\text{\eqref{eq:phi0phi1Kant0}}}{\geq}
            \phi_1(\eta(t_1)) - \phi_0(\eta(t_0))
            \overset{\text{\eqref{eq:Kantequality}}}{\geq}
            \action[\eta],
          \end{equation}
          the curve $\eta$ is a minimizing extremal.
          Since minimizing extremals solve the Euler{--}Lagrange
          equation, any two curves extending $\eta$ that
          satisfy~\eqref{eq:Kantequality}
          are restrictions of a single solution of the
          Euler{--}Lagrange equation, up to a translation of the
          domain.

          \medskip
          Let $V$ and $T$ be as in
          Assumption~\ref{item:compactness} for the compact set
          $\cX_0\times \cX_1$.
          Let $t_1'$ denote the supremum of all $s\geq t_1$ such that the curve $\eta$
          can be extended to a curve $\sigma$ satisfying~\eqref{eq:Kantequality} on the interval $[t_0,s]$ with $\sigma(s) \in \cX_1$.
          Note that $s - t_0\le T$ for every such $s$ by Assumption~\ref{item:compactness}, and in particular $t_1'$ is finite.
          Thus the curve $\eta$ can be
          extended to a curve $\sigma:[t_0,t_1')\to M$,
          which is again a solution to the Euler{--}Lagrange equation.
          By the definition of $V$, the velocities $(\dot\sigma(s))_{s\in[t_0,t_1')}$ are contained in a compact subset of $\interior\dom L$,
          and thus the curve $\sigma$ extends to a solution
          $\sigma:[t_0,t_1']\to M$ of the Euler{--}Lagrange equation.
          Since $\cX_1$ is closed, $\sigma(t_1') \in \cX_1$.
          Choose a sequence $s_n \to t_1'$ such that
          $\sigma|_{[t_0,s_n]}$ satisfies~\eqref{eq:Kantequality}.
          Since $\phi_i$ is continuous on $\cX_i$, $i=0,1$,
          \[
          \action[\sigma]
          =
          \lim_{n\to\infty}
          \action[\sigma|_{[t_0,s_n]}]
          =
          \lim_{n\to\infty}
          \big(
          \phi_1(\sigma(s_n))
          -
          \phi_0(\sigma(t_0))
          \big)
          =
          \phi_1(\sigma(t_1'))
          -
          \phi_0(\sigma(t_0))
          ,
          \]
        so $\sigma$ satisfies~\eqref{eq:Kantequality}.
        By a similar argument, the domain of $\sigma$ can be maximally extended from below; after a translation of the domain, this yields the curve $\bar\gamma:[0,\ell]\to M$.
      \end{proof}

      We can thus define transport rays and transport sets, cf.~\cite{EvGa99,FeMc02,Kl17}.
      \begin{definition}[(Upper) transport rays and endpoints]\label{def:rays}
            Let $\bar\gamma:[0,\ell] \to M$ be a maximally defined curve satisfying
            \begin{equation}\label{eq:ray-Kant}
              \bar\gamma(0)\in\cX_0,
              \quad
              \bar\gamma(\ell)\in\cX_1,
              \qquad
              \text{ and }
              \qquad
              \phi_1(\bar\gamma(\ell)) - \phi_0(\bar\gamma(0))
              \geq
              \action[\bar\gamma],
            \end{equation}
            see Lemma~\ref{lem:calibrated}.
            The curve $ \gamma : = \bar\gamma\vert_{(0,\ell)} $
            is called a \emph{transport ray}
            and the collection of transport rays is denoted by $\Gamma$.
            The length $\ell$ of the domain of a transport ray $\gamma$
            will be denoted by $\ell_\gamma$.
            A curve of the form
            $\bar\gamma\vert_{(t,\ell)}$ is called an \emph{upper
            transport ray}.
            The points $\bar\gamma(0)$ and $\bar\gamma(\ell)$
            are called the \emph{endpoints} of the transport ray $\gamma$.
          \end{definition}
          \begin{definition}[Transport sets]
          A Borel set which is a union of transport rays is called a
          \emph{transport set}.
          A Borel set which is a union of upper transport
          rays is called an \emph{upper transport set}.
        \end{definition}

        The union of all transport rays is denoted by $\cT$; thus
        \[
          \cT : =
          \bigcup_{\gamma \in \Gamma}\gamma((0,\ell_\gamma))
          .
        \]
        For every $k \ge 1$ we write
        \begin{equation}\label{eq:Tk-def}
            \cT_k : = \bigcup_{\gamma \in \Gamma}\gamma\left(\left(2^{-k},\ell_\gamma - 2^{-k}\right)\right).
        \end{equation}
        Note that
        \[
            \cT = \bigcup_{k\ge 1} \cT_k.
        \]
        
        Set
        \begin{equation}\label{eq:u-def}
          u :  M \to [-\infty,\infty]
          \qquad
          u(x) :=
          \inf_{x_0 \in \cX_0}\bigl[\phi_0(x_0) + \cc(x_0,x)\bigr]
          .
        \end{equation}
        
        A function $u : M \to [-\infty,\infty]$ is said to be \emph{$L$-dominated} if
        \begin{equation}\label{eq:udominated}
          u(x') \le u(x) + \cc(x,x')
          \qquad
          \text{for all $x,x' \in M$,}
        \end{equation}
        with the convention $+\infty-\infty = -\infty + \infty =
        +\infty$.

        \begin{definition}[cf.~\cite{Fathi08,BernardBuffoni06}]
          A piecewise-$C^1$ curve $\gamma:I\to M$ defined on an interval $I\subseteq\RR$ is said to be \emph{$(u,L)$-calibrated}
          if $u(\gamma(t))$ is finite for all $t\in I$,
          \begin{equation}
            \label{eq:u-calibrated}
            u(\gamma(t))
            -u(\gamma(s))
            =
            \action[\gamma|_{[s,t]}]
            ,
            \qquad
            \text{ for all }
            s<t \in I
            ,
          \end{equation}
          and $\dot\gamma(t)\in\interior\dom L$ for all $t\in I$.
        \end{definition}

        When $u$ and $L$ are clear from the context, we simply say that $u$ is \emph{dominated} and that $\gamma$ is \emph{calibrated}.
        Note that the restriction of a calibrated curve to any compact subinterval of its domain is a
        minimizing extremal.

        \medskip
        In the following lemma we collect some well-known facts about the function $u$ and transport rays. For most of the details we will provide the reader with references to existing proofs in the literature.
        A subset of $M$ is said to have \emph{measure zero} if it has measure zero with respect to some, hence every, measure on $M$ with a smooth positive density.
      
        \begin{lemma}\label{lem:strain}{~}
            \begin{enumerate}
                \item\label{item:u-facts}
          The function $u$ defined in~\eqref{eq:u-def} is $L$-dominated,
          and every transport ray is calibrated; moreover,
          \begin{equation}
            \label{eq:phi-u-coincide}
            u = \phi_i,
            \qquad
            \text{ on }\cX_i\cap\cT,
            \qquad
            i=0,1;
          \end{equation}
          in particular, $u$ is bounded on $(\cX_0\cup\cX_1)\cap\cT$.
                \item\label{item:transportray}
                  For every $x \in \cT$ there
                  exists a unique transport ray
                  \[
                    \gamma_x \in \Gamma
                  \]
                  passing through $x$ at a unique time, i.e., there exists a
                  unique $t_x \in (0,\ell_{\gamma_x})$ such that
                  \[
                    \gamma_x(t_x) = x.
                  \]
                \item\label{item:diff}\label{item:dotgamma}\label{item:HJzero} The function $u$ is differentiable on the set $\cT$, with
                      \begin{equation}\label{eq:dotgammatcLdu}
                          \dot\gamma(t)
                          =
                          \cL\left(\de u \vert_{\gamma(t)}\right)  \qquad \text{ for all $t \in (0,\ell_\gamma)$ and every transport ray $\gamma\in\Gamma$,}
                      \end{equation}
                      and $H(\de u)\equiv 0$ on the set
                      $\cT$.

                \item\label{item:ends} The set of endpoints of transport rays has measure zero, and
                    \begin{equation}\label{eq:strainfullmeasure}
                        \mu_0\left(M\setminus \cT\right) = \mu_1\left(M\setminus \cT\right) = 0.
                    \end{equation}
                    In particular $(u,u)$ is a Kantorovich pair.

                  \item\label{item:balance}
                    If $\cT' \subseteq \cT$ is a transport set then
                    $\mu_0(\cT') = \mu_1(\cT')$, and if $\cT^+
                    \subseteq \cT$ is an upper transport set then
                    $\mu_0(\cT^+) \le \mu_1(\cT^+)$.
                \item\label{item:borel}
                  The set $\cT$ is Borel, and the functions $x\mapsto t_x$ and $x\mapsto\ell_{\gamma_x}$ and the map $x\mapsto\dot\gamma_x(t_x)$ are Borel on $\cT$; in particular, the sets $\cT_k$ are Borel. 
                \item\label{item:lipschitz}
                  The differential $\de u : \cT \to T^*M$ is locally Lipschitz on $\cT_k$ for every $k \ge 1$.
            \end{enumerate}
        \end{lemma}
        \begin{proof}
        The function $u$ is $L$-dominated by the triangle
        inequality~\eqref{eq:triangleineq}.
        By~\eqref{eq:phi0phi1normalization}
        and~\eqref{eq:u-def}, 
        \begin{equation}\label{eq:phi-u-coincide1}
          u = \phi_1 \qquad \text{on $\cX_1$.}
        \end{equation}

        Let $\gamma \in \Gamma$ be a transport ray, and let $\bar\gamma:[0,\ell]\to M$ be as in Definition~\ref{def:rays}, so that $\gamma = \bar\gamma\vert_{(0,\ell)}$, and set $x_0 := \bar\gamma(0) \in \cX_0$ and $x_1 := \bar\gamma(\ell) \in \cX_1$.
        Then, for all $0<t<\ell$,
        \begin{align*}
          \phi_1(x_1)
          &\overset{\text{\eqref{eq:phi-u-coincide1}}}{=}
          u(x_1)
          \overset{\text{\eqref{eq:udominated}}}
          {\leq}
            u(\gamma(t)) +
            \action[\bar\gamma|_{[t,\ell]}]
            \\
            \overset{\text{\eqref{eq:u-def}}}
            &
            \leq
            \phi_0(x_0) +
            \action[\bar\gamma|_{[0,t]}]
            +
            \action[\bar\gamma|_{[t,\ell]}]
            =
            \phi_0(x_0) +
            \action[\bar\gamma]
            \overset{\text{\eqref{eq:ray-Kant}}}
            \leq
            \phi_1(x_1),
        \end{align*}
        whence equality holds throughout; in particular,
        \begin{equation}\label{eq:u-along-curve}
          u(\gamma(t)) = \phi_1(x_1) - \action[\bar\gamma|_{[t,\ell]}] < \infty,
          \qquad 0<t<\ell.
        \end{equation}
        
        Subtracting~\eqref{eq:u-along-curve} at two times $0<s<t<\ell$ shows that $\gamma$ is calibrated.
        If $\gamma(t) \in \cX_0$ for some $0<t<\ell$, then
        \[
          u(\gamma(t))
          \overset{\text{\eqref{eq:u-along-curve}}}{=}
          \phi_1(x_1) - \action[\bar\gamma|_{[t,\ell]}]
          \overset{\text{\eqref{eq:cdef}}}{\le}
          \phi_1(x_1) - \cc(\gamma(t),x_1)
          \overset{\text{\eqref{eq:phi0phi1Kant0}}}{\le}
          \phi_0(\gamma(t));
        \]
        since every point of $\cX_0\cap\cT$ is of this form, 
        it follows that $u \le \phi_0$ on $\cX_0\cap\cT$.
        Conversely, for every $x\in\cX_0$,
        \[
          \phi_0(x)
          \overset{\text{\eqref{eq:phi0phi1normalization}}}{=}
          \max_{y\in\cX_1}\left[\phi_1(y) - \cc(x,y)\right]
          \overset{\text{\eqref{eq:phi-u-coincide1}}}{=}
          \max_{y\in\cX_1}\left[u(y) - \cc(x,y)\right]
          \overset{\text{\eqref{eq:udominated}}}{\le}
          u(x),
        \]
        which proves~\eqref{eq:phi-u-coincide} for $i=0$.

          \medskip
          For~\ref{item:transportray}, see e.g.~\cite[p.~156]{Figalli07}.
          For~\ref{item:diff}, see~\cite[Th.~4.3.8]{Fathi08} or~\cite[Prop.~3.2]{Contreras01}
          (the proof does not require $L$ to be superlinear or to be finite everywhere, since it only involves variations of $\gamma_x$ with velocities near $\dot\gamma_x\in\interior\dom L$).
          For~\ref{item:ends} see~\cite[Lems.~3.8 and~4.2]{Kl17}, and for~\ref{item:balance} see~\cite[Lem.~4.4 and Cor.~4.5]{Kl17}.

          \medskip For $v\in\interior\dom L$, let $\gamma_v$ denote
          the solution of the Euler{--}Lagrange equation~\eqref{eq:EL}
          with initial velocity $v$, and let $\mathcal{R}$ denote the
          set of pairs $(v,\ell)\in\interior\dom L\times{[0,\infty)}$
          such that $\gamma_v$ is defined on ${[0,\ell]}$ and
          $\gamma_v\vert_{[0,\ell]}$
          satisfies~\eqref{eq:Kantequality}.
          The set $\mathcal{R}$ is compact by Lemma~\ref{lem:calibrated}, Assumption~\ref{item:compactness} applied to the set $\cX_0\times\cX_1$, and the continuity of $\phi_0$ and $\phi_1$.
          By Lemma~\ref{lem:calibrated} and~\ref{item:transportray}, for all $c,c'\ge0$,
          \[
            \{x\in\cT : t_x>c,\ \ell_{\gamma_x}-t_x>c'\}
            =
            \{\gamma_v(t) : (v,\ell)\in\mathcal{R},\ c<t<\ell-c'\},
          \]
          and $\dot\gamma_x(t_x)=\dot\gamma_v(t)$ whenever $x=\gamma_v(t)$ with $(v,\ell)\in\mathcal{R}$ and $0<t<\ell$.
          Hence these sets, and in particular the set $\cT$, as well as the graph of the map $x\mapsto\dot\gamma_x(t_x)$, are continuous images of $\sigma$-compact sets, and thus Borel, and~\ref{item:borel} follows.

        \medskip
        It remains to prove~\ref{item:lipschitz}, cf.~\cite[Th.~4.11.5]{Fathi08} and~\cite[Prop.~15]{FeMc02}. Let $V$ be the compact set given by Assumption~\ref{item:compactness} for $K = \cX_0\times\cX_1$.
        By Lemma~\ref{lem:calibrated}, every transport ray is the restriction of a minimizing extremal joining $\cX_0$ to $\cX_1$,
        and therefore $\dot\gamma_x(t) \in V$ for every $x \in \cT$ and $t \in (0,\ell_{\gamma_x})$.
        Let $k \ge 1$ and let $z\in \cT_k$.
        The same proof as in~\cite[Lem.~4.12]{Assouline24}
        (cf.~\cite[Th.~3.2]{Fathi03}) shows that there exist a
        coordinate neighborhood $U \ni z$ and a constant $C$ such that
          \begin{equation}\label{eq:sandwich}
            \bigl|u(x') - u(x) - \de u\vert_x(x'-x)\bigr| \le C\,|x'-x|^2
            \qquad \text{for all $x \in U\cap\cT_k$ and all $x' \in U$,}
          \end{equation}
          where 
          the chart is used to view $x,x'$ as points in $\RR^n$ (in~\cite[Lem.~4.12]{Assouline24}, local Lipschitz continuity of $u$ is used in order to confine the velocities of the curves $\{\gamma_x\vert_{[t_x - 2^{-k},\,t_x + 2^{-k}]}\}_{x \in U\cap\cT_k}$ to a compact subset of $\interior\dom L$; here $V$ is such a subset so local Lipschitz continuity is not needed).
          Given~\eqref{eq:sandwich}, the proof proceeds as in~\cite[Lem.~4.13]{Assouline24} (cf.~\cite[Lem.~3.8]{Fathi03}).
        \end{proof}

    \section{Tonelli regularization}\label{sec:tonellireg}

    The needle
    decomposition result Theorem~\ref{thm:main} was
    established for \emph{Tonelli} Lagrangians
    in~\cite[Th.~5.9]{Assouline24}.
    In this section we prove that, locally on $\cT$, we can replace
    $L$ by a Tonelli Lagrangian satisfying the standing assumption of~\cite{Assouline24}.

    \begin{proposition}\label{prop:localLagrangian}
      For every $k \ge 1$ and every $x \in \cT_k$, there exist
      a neighborhood $U^{x,k} \ni x$, a Lagrangian
      $L^{x,k} : T U^{x,k} \to \RR$,
      and a function $u^{x,k} : U^{x,k} \to \RR$
      such that the following hold:
        \begin{enumerate}
            \item\label{item:xk-standing} The Lagrangian $L^{x,k}$ is
              Tonelli and, on the manifold $U^{x,k}$, it satisfies
              the assumptions in~\cite[Sect.~2.3]{Assouline24}, namely
    \begin{itemize}
      \item
        supercriticality, i.e.,
        closed curves have positive action;
      \item
        every two points are joined by a minimizing extremal;
      \item
        for every compact set $A \subset U^{x,k}$ there exists a
        compact set $\tilde A \subset U^{x,k}$ containing every
        minimizing extremal with endpoints in $A$.
    \end{itemize}

            \item\label{item:xkL} For every $y \in \cT_k\cap U^{x,k}$, the Lagrangians $L$ and $L^{x,k}$ coincide on a neighborhood of the vector $\dot\gamma_y(t_y) = \cL \de u \vert_y$.
            \item\label{item:xk-calibrated} The function $u^{x,k}$ is $L^{x,k}$-dominated,
              and every $y \in \cT_k\cap U^{x,k}$ lies in
              the interior of a maximal $(u^{x,k},L^{x,k})$-calibrated curve
            which is a sub-arc of $\gamma_y$.
        \end{enumerate}
    \end{proposition}
    \begin{remark}
      In~\cite[Sect.~2.3]{Assouline24} it is also assumed that $0$ is a regular
      value of the energy; however, this is a consequence of
      supercriticality (cf.\ Lemma~\ref{lem:zeroreg}).
    \end{remark}

    The rest of the section is devoted to proving Proposition~\ref{prop:localLagrangian}.

    \subsection{The regularized Lagrangian}
    Fix $k \ge 1$ and $x \in \cT_k$. Recall that $\gamma_x:(0,\ell_{\gamma_x}) \to M$ is the transport ray through $x$, and that $\gamma_x(t_x)=x$ for a unique $t_x \in (2^{-k},\ell_{\gamma_x} - 2^{-k})$.
    Choose a precompact open neighborhood $\tilde U \ni x$ whose closure is contained in the domain of local coordinates in which
    \begin{align}\label{eq:tildeUdef}
      \gamma_x(t)=(t,0,\dots,0)
      \,\,
      \text{for all $t$ sufficiently close to $t_x$.}
    \end{align}

    Assumption~\ref{item:positivezerosection} implies that $\inf\{H(p):p\in T_x^*M\}=-L(\bz\vert_x)<0$.
    Since $H$ is smooth on the interior of its domain, it follows that,
    by taking $\tilde U$ sufficiently small, 
    we can find a smooth function $\zeta:M \to \RR$
    such that $\de \zeta\vert_y \in \interior\dom H$ for every $y$ in the closure of $\tilde U$, and
    \begin{equation}\label{eq:Hp0}
      h_0
      :=
      \sup_{\tilde U} \, H(\de \zeta)<0.
    \end{equation}
    
    Let $g_0$ denote the Euclidean metric on $\tilde U$ induced by the
    local coordinate chart,
    and let $g$ be a complete Riemannian metric on $\tilde U$ satisfying
    $g\geq g_0$.
    Denote by $|\cdot|_0$ and $|\cdot|_g$ the norms induced by the metrics $g_0$ and $g$, respectively.
    Then
    \begin{equation}\label{eq:gg0comparison}
      |v|_0 \le |v|_g
      ,
      \quad
      \forall v \in T\tilde U
      \qquad
      \text{ and }
      \qquad
      |p|_g \le |p|_0
      ,
      \quad
      \forall p \in T^*\tilde U
      .
    \end{equation}
    
    We regularize the Lagrangian $L$ in two steps: we first add a superlinear term to the Hamiltonian, so that the Lagrangian becomes finite everywhere, and then add a superlinear term and a uniformly convex term to the resulting Lagrangian.
    Let
    $\vphi : \RR \to {[0,\infty)}$
    be a smooth, convex, nondecreasing, superlinear function
    vanishing on the interval $[0,1]$.
    Let $\eps > 0$ be determined later and write
    $\vphi_\eps(t) : = \vphi(\eps t)$.
    Define
    \begin{equation}\label{eq:hatHdef}
      \hat H(p) : = H(p) + \vphi_\eps(|p-\de \zeta|_0), \qquad p \in T^*\tilde U
    \end{equation}
    and let $\hat L$ be the convex conjugate of $\hat H$:
    \begin{equation}\label{eq:hatLdef}
      \hat L(v)
      : =
      \sup_{p \in T_y^*\tilde U}\Bigl[p(v)
      -
      \hat H(p)\Bigr]
      ,
      \qquad
      v \in T_y\tilde U,
      y \in \tilde U.
    \end{equation}
    In particular, since $\vphi_\eps(0) = 0$,
    \begin{equation}\label{eq:hatLlower}
      \hat L(v)
      \overset{\text{\eqref{eq:hatLdef}}}{\geq}
      \de\zeta(v)
      -
      \hat H(\de\zeta\vert_y)
      \overset{\text{\eqref{eq:hatHdef}}}{=}
      \de\zeta(v)
      -
      H(\de\zeta\vert_y)
      \overset{\text{\eqref{eq:Hp0}}}{\geq}
      \de\zeta(v)
      -
      h_0
      ,
      \qquad
      v \in T_y\tilde U,
      y \in \tilde U.
    \end{equation}
    Write
    \[
      v_x
      : =
      \dot\gamma_x(t_x)
      =
      \cL(\de u\vert_x)
    \]
    and let $V$ be a smooth, compactly-supported vector field on
    $\tilde U$ satisfying 
    \[  
     V\vert_x = v_x.
    \]
    
    Let $\chi:\tilde U \to [0,1]$ be a smooth function vanishing on a neighborhood of $x$
    and equal to $1$ outside a compact subset of $\tilde U$.
    Finally, define the Lagrangian $L'$ by
    \begin{equation}
      \label{eq:lprime-definition}
      L'(v)
      :=
      \hat L(v)
      +
      \vphi_\eps\left(|v-V|_g^2\right)
      +
      \chi(y)  \,|v-V|_g^2,
      \qquad
      v \in T_y \tilde U,
      y \in \tilde U
      .
    \end{equation}

    \begin{lemma}\label{lem:LonB}
      There exist $\eps>0$ and a precompact open neighborhood $N \subseteq T\tilde U$ of $v_x$ with $\overline N \subseteq \interior\dom L$ such that $L' = L$ on $N$.
    \end{lemma}
    \begin{proof}
    Assume $\eps \le 1$, so that $\vphi_\eps$ vanishes on $[0,1]$.
    Since $\chi$ vanishes near $x$ and $V\vert_x = v_x$, there exists a precompact open neighborhood $N \subseteq T\tilde U$ of the tangent vector $v_x$, independent of $\eps$, such that
    \begin{equation}\label{eq:phivanish}
      \overline N
      \subseteq
      \interior\dom L\cap T \tilde U
      \qquad \text{and} \qquad L' = \hat L \quad \text{on the set $N$.}
    \end{equation}
    
      Write $N^* : = \cL^{-1}(N) \subseteq T^*\tilde U$.
      Since $\cL$ is a
      diffeomorphism, the set $N^*$ is a precompact open set.
      By taking $\eps < \inf_{p \in N^*}|\de \zeta - p|_0^{-1}$, we can guarantee that $\hat H=H$ on $N^*$. Let $v \in N$ and denote the basepoint of $v$ by $y$.
      Since $\vphi$ is nonnegative,
      \[
        p(v) - \hat H(p) \le p(v) - H(p) \le L(v)
        \qquad
        p \in T_{y}^*\tilde U,\,\,y \in \tilde U,
      \]
      and equality is attained in both inequalities when $p = \cL^{-1}(v) \in N^*$.
      Taking supremum over $p$ and using~\eqref{eq:hatLdef}, we see that $\hat L(v) =L(v)$ on the set $N$. Together with~\eqref{eq:phivanish}, this implies that $L'(v) = L(v)$ for all $v \in N$.
    \end{proof}
    From now on we fix $\eps>0$ and $N$ as in Lemma~\ref{lem:LonB}.
    Let $\action'$ and $\cc'$ denote the action and cost associated to the Lagrangian $L'$, respectively.

    \begin{lemma}\label{lem:tonelli}
      The Lagrangian $L'$ is smooth on $T\tilde U$, uniformly convex:
      \begin{equation}
        \label{eq:lprime-convex}
        \inf_{v \in T\tilde U}
        \inf_{|w|_g = 1}D^2L'\vert_v(w,w) > 0,
      \end{equation}
      uniformly superlinear: for every $A > 0$ there exists $C_A$ such that
      \begin{equation}
        \label{eq:lprime-superlinear}
        L'(v)
        \ge
        A|v|_g - C_A
        ,
        \qquad
        v \in T\tilde U,
      \end{equation}
      and Ma\~n\'e supercritical: every absolutely continuous
      closed curve $\gamma:[0,T] \to \tilde U$ satisfies
      \begin{equation}
        \label{eq:lprime-supercritical}
        \int_0^TL'(\dot\gamma(t))\,\de t \ge - h_0 \, T,
      \end{equation}
      where $-h_0 > 0$ is given by~\eqref{eq:Hp0}.
    \end{lemma}

    \begin{proof}
      Write $q := |\cdot - V|_g^2$, so that the last two terms in~\eqref{eq:lprime-definition} are $\vphi_\eps \circ q$ and $\chi\,q$.
      Constants in the proof are allowed to depend on anything as long as they are uniform over $v \in T\tilde U$.

      \medskip
        By our assumptions on $L$ and $\vphi$ and by classical Legendre{--}Fenchel duality~\cite[Th.~26.3]{Ro70}, the Hamiltonian $\hat H$ is convex with a positive definite vertical Hessian on the interior of its domain, essentially smooth, and superlinear.
        Hence, again by duality~\cite[Sect.~26]{Ro70}, on each fiber the Lagrangian $\hat L$ is finite and differentiable.
        For every $v \in T\tilde U$, the supremum in the definition of $\hat L(v)$ is
        attained uniquely at $\hat\cL^{-1}v$, where $\hat\cL : = D\hat H$ is the Legendre transform associated with $\hat H$.
        Since $\hat H$ is smooth on $\interior\dom H$ and
        has positive definite vertical Hessian there, the
        inverse function theorem implies that
        $\hat\cL^{-1}$ is smooth on
        $T \tilde U$.
        Hence $\hat L$, which satisfies $\hat L(v) = (\hat\cL^{-1}v)(v) - \hat H(\hat\cL^{-1}v)$
        for every $v \in T\tilde U$, is smooth on $T\tilde U$.
        Since $\vphi$, $\chi$ and $q$ are smooth as well, so is $L'$.

        \medskip
        All three terms in~\eqref{eq:lprime-definition} are fiberwise convex, and the vertical Hessians of the last two are $\vphi_\eps''(q)\,\de q \otimes \de q + 2\vphi_\eps'(q)\,g \ge 2\eps\vphi'(\eps q)\,g$ and $2\chi\,g$. Since $\vphi$ is convex and superlinear, 
        we may fix $t_0 > 0$ such that $\vphi'(t) \ge 1$ for all $t \ge t_0$.
        Then $D^2L'\vert_v \ge 2g$ wherever $\chi = 1$, and $D^2L'\vert_v \ge 2\eps\, g$ when $q(v) \ge t_0/\eps$.
        The rest of $T\tilde U$ is contained in a compact set on which the smooth positive definite form $D^2L'$ is bounded below by a positive multiple of $g$.

        \medskip
        Since $\tilde U$ is precompact, $|\de\zeta|_g$ is bounded on $\tilde U$ by~\eqref{eq:gg0comparison}, so by~\eqref{eq:hatLlower},
        \[
          \hat L(v)
          \ge
          \de \zeta(v) - h_0
          \ge
          -\Bigl(\sup_{\tilde U}|\de\zeta|_g\Bigr)|v|_g
          ,
          \qquad
          v \in T\tilde U
          .
        \]
        
        Since $\vphi$ is superlinear, there exists $b>0$ such that $\vphi(t) \ge 2t/\eps - b$ for all $t \ge 0$; hence $\vphi_\eps \circ q \ge 2q - b$.
        Since $V$ is compactly supported, $2q(v) \ge |v|_g^2 - 2\sup|V|_g^2$.
        Therefore, by~\eqref{eq:lprime-definition},
        \[
          L'(v)
          \ge
          |v|_g^2
          -
          \Bigl(\sup_{\tilde U}|\de\zeta|_g\Bigr)|v|_g
          -
          2\sup|V|_g^2
          -
          b
          ,
          \qquad
          v \in T\tilde U
          ,
        \]
        and~\eqref{eq:lprime-superlinear} follows.

        Let $\gamma:[0,T] \to \tilde U$ be an absolutely continuous closed curve;
        then $\int_0^T\de\zeta(\dot\gamma(t))\,\de t = 0$, and therefore
        \begin{align*}
          \int_0^TL'(\dot\gamma(t))\,\de t
          & \stackrel{\text{\eqref{eq:lprime-definition}}}{\ge}
          \int_0^T\hat L(\dot\gamma(t))\,\de t
            \stackrel{\text{\eqref{eq:hatLlower}}}{\ge}
            \int_0^T
            [\de \zeta(\dot\gamma(t)) -
            h_0
            ]\,\de t
            \geq
            - h_0
            T
          .
        \qedhere
        \end{align*}
    \end{proof}

    \begin{lemma}\label{lem:minimizers}
      Every two points of $\tilde U$ are joined by a minimizing $L'$-extremal.
      Moreover, for every compact set $A \subseteq \tilde U$ there exists a compact set $\tilde A \subseteq \tilde U$ containing every minimizing $L'$-extremal with endpoints in $A$.
    \end{lemma}

    \begin{proof}
      By Lemma~\ref{lem:tonelli}, $L'$ is a Tonelli and supercritical
      Lagrangian on the complete Riemannian manifold $(\tilde U, g)$.
        Existence of minimizers therefore follows from~\cite[Prop.~4.1.11]{Sorrentino15} (in~\cite{Sorrentino15} the manifold is assumed to be compact, but~\cite[Rem.~4.1.2]{Sorrentino15} clarifies that this assumption is not needed if the Lagrangian is uniformly superlinear, as in our case).

      \medskip
        Fix $C \ge \max\{C_1,1/(-h_0)\}$, where $C_1$ is given by~\eqref{eq:lprime-superlinear} with $A = 1$, and $h_0$ is as in~\eqref{eq:Hp0}.
        Let $A \subseteq \tilde U$ be a compact set; we may assume that $A$ is a closed $g$-ball of radius $R > 0$ centered at some $z \in \tilde U$.
        By Lemma~\ref{lem:tonelli}, the Lagrangian $L'$ satisfies Assumptions~\ref{item:convex}--\ref{item:domaincompat} on $\tilde U$, and it satisfies~\ref{item:positivezerosection} since $L' \ge \hat L \ge -h_0 > 0$ on the zero section by~\eqref{eq:hatLlower}.
        Hence the cost $\cc'$ is upper semicontinuous by Lemma~\ref{lem:cost-usc}, so
        \[
          D
          :=
          \max_{y,y' \in A}\cc'(y,y')
          <
          \infty.
        \]
        Fix $y,y' \in A$ and let $\gamma:[0,T] \to \tilde U$ be a minimizing $L'$-extremal joining $y$ to $y'$.
        By concatenating a minimizing $L'$-extremal from $z$ to $y$, the curve $\gamma$, and then a minimizing $L'$-extremal from $y'$ to $z$, we obtain a closed curve $\eta$ defined on an interval of length at least $T$.
        By~\eqref{eq:lprime-supercritical} and our choice of $C$,
        \[
          D
          \ge
          \cc'(y,y')
          =
          \action'[\gamma]
          =
          \action'[\eta] - \cc'(z,y) - \cc'(y',z)
          \ge
          \action'[\eta] - 2D
          \ge
          \frac{T}{C} - 2D,
        \]
        so $T \le 3CD$.
        Moreover,
        \[
          D
          \ge
          \action'[\gamma]
          \overset{\text{\eqref{eq:lprime-superlinear}}}\ge
          \int_0^T|\dot\gamma(t)|_g\,\de t - TC
          \ge
          \Len\gamma - TC
          \ge
          \Len\gamma- 3C^2D,
        \]
        so $\Len\gamma \le (1+3C^2)D$, and therefore $\gamma$ is contained in the closed $g$-ball of radius $R + (1+3C^2)D$ centered at $z$.
    \end{proof}

    \subsection{A convex neighborhood and a transverse hypersurface}

    \begin{lemma}\label{lem:convex}
      The point $x$ admits a neighborhood basis $\sW$ of open sets such that, for every $W \in \sW$, every two points in $W$ are joined by a unique minimizing $L'$-extremal
      which is contained in $W$,
      depends continuously on its endpoints,
      and is the only solution of the Euler{--}Lagrange equation for $L'$ with zero energy which joins these points within $W$ (up to a shift of the parameter).
    \end{lemma}
    \begin{proof}
      On the zeroth energy level of $L'$, the Euler{--}Lagrange flow
      coincides up to time repara\-me\-tri\-zation with the geodesic
      flow of a Finsler metric $\mathcal{F}$,
      see~\cite[Th.~1, p.~3314]{IturriagaSanchezMorgado00}.
      Hence, by Whitehead's
      theorem~\cite[Cor.~14.1.2]{shen2001differential}, there
      exists a neighborhood basis $\sW'$ of $x$ consisting of convex sets.
      Here ``convex'' means that, given two points
      $y,y'\in W\in \sW'$, there exists a unique minimizing geodesic joining them
      and contained in $W$, and such a geodesic is the unique solution
      in $W$ of the geodesic equation with $y$ and $y'$ as boundary
      data.
      Fix $W_0\in \sW'$.
      By~\cite[Lem.~2.9]{Assouline24}, whose assumptions hold by Lemmas~\ref{lem:tonelli} and~\ref{lem:minimizers}, there exists a neighborhood
      $W_1 \subseteq W_0$ of $x$ such that for every $y,y' \in W_1$,
      there exists a unique minimizing extremal $\sigma$ joining $y$ to $y'$ in $\tilde U$, and $\sigma$ is contained in $W_0$.
      Take 
      \[
        \sW:=\{W\in\sW'\colon W\subset W_1\}.
      \]
      
      We check that $\sW$ has the desired properties.
      Fix $y,y'\in W\in \sW$.
      Let $\sigma$ be the unique minimizing extremal (for $L'$) joining
      $y$ to $y'$.
      Then $\sigma$ is contained in $W_0$ and solves the Euler{--}Lagrange
      equations, and therefore, up to a reparametrization, solves the
      geodesic equation of the Finsler metric $\mathcal{F}$.
      By uniqueness (in $W_0$) of the solutions of the geodesic
      equations, $\sigma$ is a (not necessarily constant-speed) geodesic.
      Since $y,y'\in W\in \sW'$, $\sigma$ is contained in $W$.
      If $\sigma'$ is another solution of the Euler{--}Lagrange equation for $L'$ with zero energy joining $y$ to $y'$ within $W$, then, up to reparametrization, $\sigma'$ is also a geodesic of $\mathcal{F}$ joining $y$ to $y'$ within $W$, and hence has the same image as $\sigma$; since the zero energy level meets each ray in $T_{y}\tilde U$ in at most one point, the parametrization is uniquely determined up to a shift, so $\sigma' = \sigma$ up to a shift of the parameter.

      \medskip
      Continuous dependence on the endpoints follows from the
      inverse function theorem applied to the map 
      \[
        \Psi:\{E'=0\}\times(0,r)\to (\tilde U\times \tilde U)\setminus\Delta
        \qquad
        \Psi(v,t): = (\pi(v),\gamma_v(t)),
      \]
      for $r > 0$ small enough, where $E'$ denotes the energy of the Lagrangian $L'$, $\Delta\subseteq M\times M$ denotes the diagonal, and $\gamma_v$ denotes the solution to Euler-Lagrange with
      initial velocity $v$.
    \end{proof}
    The local function $u^{x,k}$ in Proposition~\ref{prop:localLagrangian} will be constructed using
    a Hopf-Lax formula for the Lagrangian $L'$, see~\eqref{eq:uprimedef} below where it is denoted by $u'$. 
    In the following lemma we construct a hypersurface which will carry the initial data for $u'$,
    and establish some facts which will guarantee that $u$ and $u'$ share the same transport rays.

    \begin{lemma}\label{lem:localU}
      There exist a neighborhood $U\subseteq \tilde U$ of the point
      $x$, a smooth compact hypersurface with boundary
      $\Sigma\subseteq \tilde U$ and $0 < \delta < 2^{-k-1}$ such
      that:
    \begin{enumerate}[label=\textup{(U\arabic*)},ref=\textup{(U\arabic*)}]
    \item\label{item:U-convex} on the manifold $U$, the Lagrangian
      $L'$ satisfies the assumptions of~\cite[Sect.~2.3]{Assouline24}
      (see
      Proposition~\ref{prop:localLagrangian}\ref{item:xk-standing});
      \item\label{item:U-window} every $y \in U\cap\cT_k$ satisfies $\gamma_y(t) \in \tilde U$ and $\dot\gamma_y(t) \in N$ for all $t \in (t_y-\delta,\,t_y+\delta)$, and $\gamma_y(t_y \pm \delta/2) \notin U$;
      \item\label{item:U-Lprimevel} for every
        $z\in\Sigma$ and $y\in U$, there exists a
        unique minimizing $L'$-extremal (in $\tilde U$) joining $z$ to $y$ and its
        velocity lies in $N$;
      \item\label{item:costs-comparison}
      $\cc(z,y)\leq \cc'(z,y)$ for all $z\in \Sigma$ and $y\in U$;
      \item\label{item:U-cross} for every $y \in U\cap\cT_k$ there exists $s_y \in (t_y - \delta,\,t_y)$ such that $\gamma_y(s_y) \in \Sigma\cap\cT_{k+1}$ and $\cc(\gamma_y(s_y),y) = \cc'(\gamma_y(s_y),y)$ (see Figure~\ref{fig:locallagrangian});
      \item\label{item:U-ubounded}
      the function $u$ is bounded on $\Sigma$.
    \end{enumerate}
    \end{lemma}
    \begin{proof}
      Since $N$ is a precompact neighborhood of $v_x = \dot\gamma_x(t_x) = \cL(\de u\vert_x) \in T\tilde U \setminus\bz$,
      by local Lipschitz continuity of $\de u$ on $\cT_k$ and
      smoothness of the Euler{--}Lagrange flow of $L$ and $L'$,
      the first part of property~\ref{item:U-window} can be guaranteed by taking $U$ and $\delta$ sufficiently small.
      Since the transport ray $\gamma_x$ passes through $x$ at a unique time, $\gamma_x(t_x\pm\delta/2) \neq x$; since the map $y \mapsto \gamma_y(t_y\pm\delta/2)$ is continuous on $\cT_k$ by Lemma~\ref{lem:strain}\ref{item:dotgamma} and the continuity of $\de u$ on $\cT_k$, the second part of~\ref{item:U-window} can be guaranteed by shrinking $U$, while keeping $\delta$ fixed.

      \medskip
      Let $W \in\sW$ be a neighborhood from the basis provided by
      Lemma~\ref{lem:convex}.
      Let
      \begin{equation}\label{eq:sxdef}
        s_x \in (t_x - \delta/4,t_x)
      \end{equation}
      and let $\Sigma$ be an $(n-1)$-dimensional hypersurface
      passing through $\gamma_x(s_x)$ and transversal to
      $\dot\gamma_x(s_x)$
      (see Figure~\ref{fig:locallagrangian}).
      Such a hypersurface can be constructed in the local
      coordinates we fixed on $\tilde U$, see~\eqref{eq:tildeUdef}, by taking
      \[
        \Sigma=\{s_x\}
        \times
        \overline{B_{\delta'}(0)}
        ,
      \]
      for some small enough $\delta'>0$.
      By taking $s_x$ sufficiently close to $t_x$ and $\delta'$
      sufficiently small, and shrinking $U$ if necessary,
      we can guarantee that the curve $\gamma_x\vert_{[s_x,t_x]}$,
      the hypersurface $\Sigma$, and the neighborhood $U$
      are all contained in the neighborhood $W$.
      By Lemma~\ref{lem:strain}\ref{item:diff}, the function $u$
      is differentiable at the point
      $\gamma_x(s_x) \in \cT$, hence $u$ is bounded above on a
      neighborhood of $\gamma_x(s_x)$; by taking $\delta'$ smaller if necessary, we can guarantee that $\Sigma$ is contained in
      this neighborhood.
      Since $\gamma_x$ is a minimizing $L$-extremal, it has zero energy.
      Since moreover $\dot\gamma_x(t) \in N$ for $t \in [s_x,t_x]$ by~\ref{item:U-window}, Lemma~\ref{lem:LonB} implies that the segment
      $\gamma_x\vert_{[s_x,t_x]}$ is also a solution to the Euler{--}Lagrange
      equations for the Lagrangian $L'$ with zero $L'$-energy.
      Thus, by the uniqueness clause of Lemma~\ref{lem:convex}, the curve $\gamma_x\vert_{[s_x,t_x]}$ is
      a minimizing extremal for the Lagrangian $L'$ as well,
      and by the continuity clause of Lemma~\ref{lem:convex}, if $\delta'$ and $U$ are chosen sufficiently small,
      then $\dot\gamma\in N$ for every minimizing $L'$-extremal
      $\gamma$ joining $\Sigma$ to $U$, i.e.,
      property~\ref{item:U-Lprimevel} holds.

      \medskip
      If $z \in \Sigma$ and $y \in U$, and $\sigma$ is the minimizing $L'$-extremal joining $z$ to $y$, then $\dot\sigma \in N$ by~\ref{item:U-Lprimevel}, and Lemma~\ref{lem:LonB} gives
      \begin{equation}
        \label{eq:costs-comparison}
        \cc(z,y)
        \leq
        \action[\sigma]
        =
        \action'[\sigma]
        =
        \cc'(z,y)
      \end{equation}
      i.e.,\ property~\ref{item:costs-comparison} holds.

      \medskip
      By the local Lipschitz continuity of $\de u$ 
      on $\cT_k$, continuous dependence of the solution of the
      Euler{--}Lagrange equation,
      transversality of
      $\Sigma$ to $\dot\gamma_x(s_x)$, and~\eqref{eq:sxdef},
      by shrinking $U$ further we can guarantee that for all
      $y\in U\cap  \cT_k$, the ray $\gamma_y$
      intersects $\Sigma$ at a time $s_y$ satisfying
      \begin{equation}
        t_y-s_y\in (\tfrac{t_x-s_x}{2},\tfrac{\delta}{2})
        \subseteq (0,2^{-k-2})
        .
      \end{equation}
      Since $t_y-s_y>0$, the ray
      $\gamma_y$ intersects $\Sigma$ before reaching $U$.
      On the other hand, since $y\in \cT_k$ and
      $s_y\geq t_y-2^{-k-2}$, we have that $\gamma_y(s_y)\in \cT_{k+1}$ (recall~\eqref{eq:Tk-def}).
      Since $s_y \in (t_y-\delta,\,t_y)$, the velocity $\dot\gamma_y$ lies in $N$ on $[s_y,t_y]$ by~\ref{item:U-window}, and Lemma~\ref{lem:LonB} gives
      \begin{equation*}
        \cc\left(\gamma_y(s_y),y\right)
        =
        \action\left[\gamma_y\vert_{[s_y,t_y]}\right]
        =
        \action'\left[\gamma_y\vert_{[s_y,t_y]}\right]
        \geq\cc'\left(\gamma_y(s_y),y\right)
        \overset{\text{\eqref{eq:costs-comparison}}}{\geq}
        \cc(\gamma_y(s_y),y)
        ,
      \end{equation*}
      proving~\ref{item:U-cross}.

      \medskip
      It remains to prove~\ref{item:U-ubounded}.
      By the choice of $\delta'$, the function $u$ is bounded above on $\Sigma$.
      For the lower bound, let $z \in \Sigma$. By~\ref{item:costs-comparison} and the fact that $u$ is $L$-dominated, 
      \[
        u(x)
        \le
        \cc(z,x) + u(z)
        \le
        \cc'(z,x) + u(z).
      \]
      Since $u(x)$ is finite, and since the function $\cc'(\cdot,x)$ is finite and continuous on $\tilde U$ and therefore bounded above on the compact set $\Sigma$, it follows that $u \ge u(x) - \max_{\Sigma}\cc'(\cdot,x)$ on $\Sigma$.
      Finally, we replace $U$ by some $U_0 \in \sW$ with $U_0 \subseteq U$; this preserves~\ref{item:U-window}--\ref{item:U-ubounded} and gives~\ref{item:U-convex}.
      Indeed, existence of minimizers is trivial.
      Supercriticality is given by Lemma~\ref{lem:tonelli} and it
      passes trivially to subsets.
      Fix a compact set $A\subset U_0$; the union of all minimizing
      $L'$-extremals with endpoints in $A$ is compact, since it is the
      image of a compact set under a continuous map.
    \end{proof}

    \subsection{Proof of Proposition~\ref{prop:localLagrangian}}

    \begin{figure}
        \centering
        \begin{tikzpicture}[scale=1.3, dot/.style={circle, fill, inner sep=0pt, minimum size=3.5pt}]

            \foreach \dx in {-0.4,-0.3,-0.2,-0.1,0,0.1,0.2,0.3,0.4} {
                    \draw[blue, thick, opacity=0.2]
                    plot[smooth]
                    coordinates {
                            (.8*\dx,1.2)
                            (.9*\dx+0.06,1.8)
                            (\dx,2.5)
                            (\dx-0.06,3.1)
                            (1.2*\dx,3.8)
                        };
                }
            \node[blue, opacity=0.3] at (.7,4) {$\cT_k$};

            \draw[blue!60, thick]
            plot[smooth]
            coordinates {
                    (0,1.2)
                    (0.06,1.8)
                    (0,2.5)
                    (-0.06,3.1)
                    (0,3.8)
                };
            \node[blue!60, above] at (0,3.8) {$\gamma_x$};

            \draw[thick, dashed] (0,2.6) circle (2.6);
            \node[below left] at (2,4) {$\tilde U$};

            \draw[thick, dashed]
            plot[smooth cycle, tension=0.65]
            coordinates {
                    (.8,2.5)
                    (0.4,3)
                    (-0.6,3)
                    (-0.8,2.4)
                    (-0.6,2)
                    (0.0,1.9)
                    (.7,2)
                };
            \node at (-1.2,2.5) {$U$};

            \draw[thick]
            plot[smooth, tension=0.9]
            coordinates {
                    (-.8,1.7)
                    (-0.4,1.6)
                    (0.4,1.7)
                    (.8,1.6)
                };
            \node[right] at (.9,1.6) {$\Sigma$};

            \node[dot] at (0,2.5) {};
            \node[right] at (0.1,2.5) {$x$};

            \node[dot] at (0.03,1.65) {};
            \node[below right] at (0,1.6) {$\gamma_x(s_x)$};

            \node[dot] at (-0.4,2.5) {};
            \node[left] at (-0.45,2.5) {$y$};

            \node[dot] at (-0.31,1.61) {};
            \node[below left] at (-0.33,1.58) {$\gamma_y(s_y)$};
        \end{tikzpicture}
        \caption{Lemma~\ref{lem:localU}}\label{fig:locallagrangian}
      \end{figure}

      We can now finish the proof of Proposition~\ref{prop:localLagrangian}.
      Fix $\Sigma$, $U$ and $\delta$ as in Lemma~\ref{lem:localU}, and take 
      \[
        U^{x,k} : = U\qquad \text{and} \qquad L^{x,k} : = L'\vert_{TU}
.
      \]
      Assertion~\ref{item:xk-standing} follows from~\ref{item:U-convex},
      and Assertion~\ref{item:xkL} follows from Lemma~\ref{lem:LonB} and~\ref{item:U-window}.

      \medskip
      We now prove assertion~\ref{item:xk-calibrated}.
      Define $u':\tilde U\to [-\infty,+\infty)$ by
      \begin{equation}\label{eq:uprimedef}
        u'(y)
        :=
        \inf_{z\in\Sigma\cap\cT_{k+1}}
        \bigl[u(z)+\cc'(z,y)\bigr]
        ,
        \qquad 
        y \in \tilde U.
      \end{equation}
      
      For each $y  \in \tilde U$, the function $\cc'(\cdot,y)$
      is bounded below on the precompact set $\Sigma\cap\cT_{k+1}$
      by continuity, and the function $u$ is bounded on $\Sigma\cap\cT_{k+1}$
      by~\ref{item:U-ubounded}, so $u'$ is real-valued.
      Note that since $\cc'(z,z)=0$, 
      \begin{equation}
        \label{eq:uprime-smaller-u-in-sigma}
        u'(z)
        \leq
        u(z)
        ,
        \qquad
        \forall
        z\in\Sigma\cap\cT_{k+1}
        .
      \end{equation}
      
      By the triangle inequality, the function $u'$ is $L'$-dominated:
      \begin{equation}\label{eq:uprimedom}
        u'(y') - u'(y) \le \cc'(y,y')
        ,
        \qquad
        \forall y,y' \in \tilde U
        .
      \end{equation}
      
      Hence the function 
      \[
        u^{x,k} : = u'\vert_U.
      \]
      is $L^{x,k}$-dominated as well, since $\cc' \le \cc^{x,k}$ on $U \times U$, where $\cc^{x,k}$ denotes the cost associated to $L^{x,k}$ (in fact, equality holds by Lemma~\ref{lem:convex}, since $U \in \sW$).

      \medskip
      Fix $y\in U\cap \cT_{k}$ and let $\bar z := \gamma_y(s_y) \in \Sigma$ with $s_y \in (t_y-\delta,\,t_y)$ as provided by~\ref{item:U-cross}.
      Since $u$ is $L$-dominated,
      \begin{align*}
        u(z)
        +
        \cc( z, y)
        \geq
        u(y)
        ,
        \qquad
        z \in \Sigma\cap\cT_{k+1},
      \end{align*}
      with equality when $z = \bar z$, because $\bar z$ lies on the
      transport ray $\gamma_y$.
      Hence
      \begin{equation}\label{eq:barzattainment}
        \min_{z\in \Sigma \cap \cT_{k+1}}
        u(z) +
        \cc( z, y)
        = 
        u(y)
        =
        u(\bar z)
        +
        \cc(\bar z, y),
      \end{equation}
      and therefore
      \begin{align*}
        u'(y)
        \overset{\text{\eqref{eq:uprimedom}}}&{\leq}
        u'(\bar z)
        +
          \cc'(\bar z, y)
          \overset{\text{\eqref{eq:uprime-smaller-u-in-sigma}}}
          \leq
        u(\bar z)
        +
        \cc'(\bar z, y)
        \stackrel{\text{\ref{item:U-cross}}}{=}
        u(\bar z)
        +
        \cc(\bar z, y)
        \\
        \overset{\text{\eqref{eq:barzattainment}}}&{=}
        \min_{z\in \Sigma \cap \cT_{k+1}}
        u(z)
        +
        \cc( z, y)
        \stackrel{\text{\ref{item:costs-comparison}}}{\leq}
        \inf_{z\in \Sigma \cap \cT_{k+1}}
        u(z)
        +
        \cc'( z, y)
        \stackrel{\text{\eqref{eq:uprimedef}}}{=}
        u'(y)
        .
      \end{align*}
      
      Thus for every $y \in U\cap\cT_k$ the infimum defining $u'(y)$ is attained at the crossing point $\bar z$ of the transport ray $\gamma_y$ with $\Sigma\cap\cT_{k+1}$, and
      \begin{equation}\label{eq:uprimeattained}
        u'(y) = u(y) = u(\bar z) + \cc'(\bar z,y)
        .
      \end{equation}

      Now let $y = \gamma(t)$ and $y' = \gamma(t')$ be points of $U \cap \cT_k$ lying on a common transport ray $\gamma$ with $0 \le t' - t < \delta$.
      Since $t' - t < \delta$, it holds that $\dot\gamma(r) \in N$,
      for all $r\in(t,t')$, by~\ref{item:U-window}, so that $\action[\gamma|_{[t,t']}] = \action'[\gamma|_{[t,t']}]$
      by Lemma~\ref{lem:LonB}.
      Since $\gamma$ is calibrated by Lemma~\ref{lem:strain}\ref{item:u-facts}, it follows that
      \[
        u'(y') - u'(y)
        \overset{\text{\eqref{eq:uprimeattained}}}{=}
        u(y') - u(y)
        \overset{\text{\eqref{eq:u-calibrated}}}
        \geq
        \action[\gamma|_{[t,t']}]
        =
        \action'[\gamma|_{[t,t']}]
        \ge
        \cc'(y,y')
        \overset{\text{\eqref{eq:uprimedom}}}{\ge}
        u'(y') - u'(y)
        ,
      \]
      so equality holds throughout; in particular
      \begin{equation}\label{eq:uprimecalibrated}
        u'(y')-u'(y) = \action'[\gamma|_{[t,t']}] = \cc'(y,y')
        .
      \end{equation}

      \medskip
      Now fix $y \in \cT_k\cap U$.
      Since $U$ is open and $t_y \in (2^{-k},\ell_{\gamma_y}-2^{-k})$, there exists $0 < \eps < \delta/2$ such that $\gamma_y((t_y-\eps,\,t_y+\eps))\subseteq U \cap \cT_k$.
      The curve $\gamma_y\vert_{(t_y-\eps,\,t_y+\eps)}$ is $(u',L')$-calibrated by~\eqref{eq:uprimecalibrated}, and hence $(u^{x,k},L^{x,k})$-calibrated.

      \medskip
      The curve $\gamma_y\vert_{(t_y-\eps,\,t_y+\eps)}$ is contained in a maximal $(u^{x,k},L^{x,k})$-calibrated curve $\eta : I \to U$
      which is a minimizing $L^{x,k}$-extremal, see~\cite[Prop.~4.3(e)]{Assouline24} or the proof of Lemma~\ref{lem:calibrated}.
      Since $L^{x,k} = L'\vert_{TU}$, it follows that $\eta$ is a minimizing $L'$-extremal.
      Since $\dot\gamma_y(t) \in N$ for all $t \in (t_y-\delta,\,t_y+\delta)$ by~\ref{item:U-window}, and $L' = L$ on $N$ by Lemma~\ref{lem:LonB}, the curve $\gamma_y\vert_{(t_y-\delta,\,t_y+\delta)}$ is an $L'$-extremal.
      The $L'$-extremals $\eta$ and $\gamma_y$ coincide on $(t_y-\eps,\,t_y+\eps)$,
      so since $\gamma_y(t_y \pm \delta/2) \notin U$ while $\eta$ is contained in $U$, 
      it follows from uniqueness of solutions of the Euler{--}Lagrange equation that $\eta$ is a sub-arc of $\gamma_y$.
      \qed

    \section{Needle decomposition}\label{sec:needle}
    
    In this section we prove our needle decomposition result, which is 
    modeled after the results of~\cite{Kl17,CaMo17,Assouline24}.
    Let $M$ be a smooth manifold, let $L:TM\to\RR\cup\{+\infty\}$ be a Lagrangian satisfying assumptions~\ref{item:convex}{--}\ref{item:positivezerosection}, let $\cA \subseteq (M\times M)\setminus\Delta$ satisfy conditions~\ref{item:attained} and~\ref{item:compactness}, and let $\mu_0$ and $\mu_1$ be absolutely continuous Borel probability measures on $M$ with compact supports $\cX_0$ and $\cX_1$ such that $\cX_0\times\cX_1\subseteq\cA$.
    Let the function $u$, 
    the set of transport rays $\Gamma$, 
    the union of transport rays $\cT$, 
    and the sets $\cT_k\subseteq\cT$
    be as in Section~\ref{sec:OT}.
    Recall that the domain of a transport ray $\gamma$ is $(0,\ell_\gamma)$,
    and that for $x\in\cT$, we denote the unique transport ray through $x$ by
    $\gamma_x$, so that $\gamma_x(t_x)=x$.

    \begin{theorem}\label{thm:main}
    Let $\mm$ be a measure on $M$ with a smooth positive density, and write $\mu_0=f_0\mm$ and $\mu_1=f_1\mm$.
    Then there exist a probability measure $\nu$ on $\Gamma$ and a family $(\mm_\gamma)_{\gamma \in \Gamma}$ of Borel measures on $M$ such that the following hold:
        \begin{enumerate}
            \item\label{item:disint}
              For every Borel set $A\subseteq M$, the function $\gamma\mapsto\mm_\gamma(A)$ is measurable, and
                  \begin{equation}\label{eq:disint}
                    \mm(A\cap\cT)
                    =
                    \int_\Gamma\mm_\gamma(A)
                    \,\nu(\de\gamma)
                    .
                  \end{equation}
            \item\label{item:needle}
                For $\nu$-a.e.\ $\gamma \in \Gamma$, the measure $\mm_\gamma$ takes the form
                \begin{equation}
                  \label{eq:m-is-a-push-forward}
                  \mm_\gamma = \gamma_{\pushf}(m_\gamma),
                \end{equation}
                where $m_\gamma$ is a measure on the domain
                $(0,\ell_\gamma)$ of $\gamma$ with a smooth, positive density
                $ e^{-\psi_\gamma} $ with respect to the Lebesgue
                measure.
            \item\label{item:massbalance}
                For $\nu$-a.e.\ $\gamma \in \Gamma$,
                \begin{equation}
                  \label{eq:mass-balance-main-theorem}
                  \int f_0 \,\de \mm_\gamma = \int f_1\, \de \mm_\gamma,
                \end{equation}
                and
                \begin{equation}
                  \label{eq:monotonicity-main-theorem}
                  \int_t^{\ell_\gamma} f_0\circ\gamma \, \de m_\gamma
                  \le
                  \int_t^{\ell_\gamma} f_1\circ\gamma\, \de m_\gamma,
                  \qquad
                  \text{ for all }
                  t \in (0,\ell_\gamma)
                  .
                \end{equation}
            \item\label{item:HJneedle}
                For $\nu$-a.e.\ $\gamma \in \Gamma$, for every
                $t \in (0,\ell_\gamma)$ there exists a smooth solution
                $u_{\gamma,t}$ to the Hamilton{--}Jacobi equation $H(\de u_{\gamma,t}) = 0$ defined on a neighborhood of $\gamma(t)$ such that
                \begin{equation}\label{eq:HJneedle}
                    \dot\gamma = (\cL \, \de u_{\gamma,t}) \circ \gamma
                    \quad \text{ and } \quad -\dot\psi_\gamma = \left(\div_\mm\left(\cL\,\de u_{\gamma,t}\right) \right)\circ\gamma
                \end{equation}
                in a neighborhood of $t$.
              \end{enumerate}
    \end{theorem}

    \begin{remark}\label{rem:curvature}
      Equation~\eqref{eq:HJneedle} immediately implies that the
      log-density of the disintegrated measures satisfies a
      Riccati-like inequality.
      Indeed, let $N\in(-\infty,\infty]\setminus[1,n)$, where $n = \dim M$; differentiating the second equation in~\eqref{eq:HJneedle} and using the first gives
      \begin{align*}
        \ddot\psi_\gamma
        -
        \frac{
        (\dot\psi_\gamma)^2
        }{N-1}
        =
        -
        (
        \de
        (\div_{\mm}
        (\cL \, \de u_{\gamma,t})
        ))(\dot\gamma)
        -
        \frac{
        (\div_{\mm}
        (\cL \, \de u_{\gamma,t})
        )^2
        }{N-1}
        \geq
        \Ric_{\mm,N}(\dot\gamma)
        ,
      \end{align*}
      where $\Ric_{\mm,N}$ is the weighted Ricci curvature for the
      Lagrangian $L$ as defined in~\cite[Sect.~3]{Assouline24},
      and the inequality follows from the Bochner inequality~\cite[Eq.~(3.15)]{Assouline24}.
      Thus, if the Ricci curvature of $L$ in the sense of~\cite{Assouline24} is bounded below, then each measure $\mm_\gamma$ satisfies a curvature-dimension condition.
    \end{remark}

    \subsection{Disintegrations}

    We define a notion of disintegration that will be used repeatedly in the proof of Theorem~\ref{thm:main}.
    We refer to, e.g.,~\cite[Appendix~A]{BianchiniCaravenna09}
    for a more detailed treatment of disintegrations.

    \begin{definition}\label{def:disint}
      Let $(X,\mm)$ be a $\sigma$-finite measure space,
      let $Y$ be a measurable space and let $Q:X\to Y$ be a measurable
      map.
      Let $\nu$ be a $\sigma$-finite measure on $Y$
      and let $(\mm_y)_{y\in Y}$ be a family of measures
      on $X$.
      We say that the pair $(\nu, (\mm_y)_{y\in Y})$ is a
      \emph{disintegration} of $\mm$ with respect to $Q$ if
      the following two conditions hold:
      \begin{enumerate}
        \item
        For $\nu$-a.e.\ $y\in Y$,
        \begin{equation}
          \label{eq:strong-consistency}
          \text{the measure }
          \mm_y
          \text{ is concentrated on }
          Q^{-1}(y).
        \end{equation}
        \item
        For every measurable subset $A\subseteq X$, the map
        $y\mapsto \mm_y(A)$ is $\nu$-measurable, and the
        disintegration formula holds:
        \begin{equation}\label{eq:def-disintegration-formula}
          \mm(A)
          =
          \int_Y\mm_y(A)\,\nu(\de y)
          .
        \end{equation}
      \end{enumerate}
      Here and in what follows, given a measure $\nu$, we say that a set or a function is \emph{$\nu$-measurable} if it is measurable with respect to the completion of $\nu$.
    \end{definition}
    \begin{remark}\label{rem:stronglyconsistent}
      What we call a disintegration is sometimes called a
        \emph{strongly consistent} disintegration (see,
        e.g.,~\cite[Def.~3.2.7]{Cavalletti17}).
    \end{remark}

    Let us collect some facts about disintegrations that are used in the proof of Theorem~\ref{thm:main}.

    \begin{proposition}[{Existence of disintegrations~\cite[Th.~5.3.1]{AmbrosioGigliSavare08}}]\label{prop:existence-disintegration}
      Let $X$ and $Y$ be Radon separable metric spaces, equipped with their Borel $\sigma$-algebras, let $\mm$ be a $\sigma$-finite measure on $X$, let $Q:X\to Y$ be a measurable map, and let $\nu$ be a $\sigma$-finite measure on $Y$ such that $Q_\pushf\mm\ll\nu$.
      Then there exists a family $(\mm_y)_{y\in Y}$ of measures on $X$ such that $(\nu,(\mm_y)_{y\in Y})$ is a disintegration of $\mm$ with respect to $Q$ and the function $y\mapsto\mm_y(A)$ is Borel for every Borel set $A\subseteq X$.
      If moreover $\nu\ll Q_\pushf\mm$, then $\mm_y\neq0$ for $\nu$-a.e.\ $y\in Y$.
    \end{proposition}

    \begin{proposition}[{Uniqueness of disintegrations~\cite[Def.~A.6 and Th.~A.7]{BianchiniCaravenna09}}]\label{prop:disint-uniqueness}
      Let $(X,\mm)$ be a $\sigma$-finite measure space whose $\sigma$-algebra is countably generated, let $Y$ be a measurable space, let $Q:X\to Y$ be a measurable map and let $\nu$ be a $\sigma$-finite measure on $Y$.
      If $(\nu,(\mm_y)_{y\in Y})$ and $(\nu,(\mm'_y)_{y\in Y})$ are disintegrations of $\mm$ with respect to $Q$, then $\mm_y=\mm'_y$ for $\nu$-a.e.\ $y\in Y$.
    \end{proposition}

    The following lemma, regarding the behavior of disintegrations under composition, follows from Definition~\ref{def:disint} by a routine argument.

    \begin{lemma}\label{lem:composite-disintegration}
      Let $X$ be a separable metric space equipped with its Borel $\sigma$-algebra, let $\mm$ be a $\sigma$-finite measure on $X$, let $Y$ and $Z$ be measurable spaces, and let $Q:X\to Y$ and $R:Y\to Z$ be measurable maps.
      Let $(\nu',(\mm_y)_{y\in Y})$ be a disintegration of $\mm$ with respect to $Q$, and let $(\nu,(\nu'_z)_{z\in Z})$ be a disintegration of $\nu'$ with respect to $R$.
      Then $\big(\nu,(\int_Y\mm_y\,\nu'_z(\de y))_{z\in Z}\big)$ is a disintegration of $\mm$ with respect to $R\circ Q$.
    \end{lemma}

    \subsection{Proof of Theorem~\ref{thm:main}}

      The proof proceeds as follows.
      We first apply Proposition~\ref{prop:existence-disintegration} to obtain a disintegration of the measure $\mm\llcorner_\cT$
      with respect to the map $x\mapsto\gamma_x$.
      We then fix $k \ge 1$
      and apply Proposition~\ref{prop:localLagrangian}
      to obtain local Tonelli Lagrangians and dominated functions at each point of the set $\cT_k$;
      the disintegration result in~\cite{Assouline24} then provides
      local disintegrations
      satisfying properties~\ref{item:needle}--\ref{item:HJneedle}
      with respect to the local Lagrangians,
      which coincide with the Lagrangian $L$ near the transport direction.
      Finally, we relate the local disintegrations to the global one
      using Proposition~\ref{prop:disint-uniqueness} and Lemma~\ref{lem:composite-disintegration}.

      \medskip
      We identify the set of transport rays $\Gamma$ with the set $\{x\in\cT:2t_x=\ell_{\gamma_x}\}$ of their midpoints, which is Borel by Lemma~\ref{lem:strain}\ref{item:borel}.
      The map
      \[
        Q:\cT\to\Gamma, \qquad Q(x):=\gamma_x
      \]
      is then Borel, since the midpoint of $\gamma_x$ depends measurably on $x$ by Lemma~\ref{lem:strain}\ref{item:transportray} and~\ref{item:borel} and the continuity of the Euler{--}Lagrange flow.
      Let $\nu$ be any Borel probability measure on $\Gamma$ that is mutually absolutely continuous with the measure $Q_\pushf(\mm\llcorner_\cT)$.
      Proposition~\ref{prop:existence-disintegration} provides a disintegration
      \begin{equation}\label{eq:disintglobal}
        (\nu,(\mm_\gamma)_{\gamma\in\Gamma})
      \end{equation}
      of $\mm\llcorner_\cT$ with respect to the map $Q$, with $\mm_\gamma\neq0$ for $\nu$-a.e.\ $\gamma\in\Gamma$.

      \medskip
      Let $k\ge1$, and for every $x\in \cT_k$ let $U^{x,k}$ be the neighborhood
      provided by Proposition~\ref{prop:localLagrangian}.
      Extract a countable subcover $(U^{i,k})_{i\in \NN}$ from the cover $(U^{x,k})_{x \in \cT_k}$,
      and denote by $L^{i,k}$ and $u^{i,k}$ the corresponding local Lagrangian
      and function furnished by
      Proposition~\ref{prop:localLagrangian}.
      The relative interior of a maximal $(u^{i,k},L^{i,k})$-calibrated curve in $U^{i,k}$ 
      will be called an \emph{$(i,k)$-ray},
      and the union of all $(i,k)$-rays will
      be denoted by $\cT^{i,k} \subseteq U^{i,k}$.

      \medskip
      Fix $i\ge1$.
      By Proposition~\ref{prop:localLagrangian},
      the Tonelli Lagrangian $L^{i,k}$ on the weighted manifold $(U^{i,k},\mm\llcorner_{U^{i,k}})$ and the $L^{i,k}$-dominated function
      $u^{i,k}$
      satisfy the assumptions of~\cite{Assouline24},
      and $\cT^{i,k}$ is the strain set of $u^{i,k}$ in the sense of~\cite[Sect.~4.1]{Assouline24}.
      Hence, by~\cite[Th.~5.9]{Assouline24}, there exists a disintegration
      \begin{equation}\label{eq:disintlocal}
        (\nu^{i,k},(\mm_{\eta}^{i,k})_{\eta})
      \end{equation}
      of the measure
      $\mm\llcorner_{\cT^{i,k}}$ with respect to
      the map sending each point of $\cT^{i,k}$ to the unique $(i,k)$-ray through it, where $\nu^{i,k}$ is a $\sigma$-finite measure on the set of $(i,k)$-rays, 
      such that the following hold for $\nu^{i,k}$-a.e.\ $(i,k)$-ray $\eta$:

      \begin{itemize}
        \item 
        The measure $\mm_\eta^{i,k}$ takes the form
        \begin{equation}
          \label{eq:local-needle}
          \mm_\eta^{i,k} =
          \eta_{\pushf}\big(e^{-\psi^{i,k}_\eta}\Leb^1\big),
        \end{equation}
        where $\psi^{i,k}_\eta$ is a smooth function on the domain of $\eta$;
        here and below $\Leb^1$ denotes the $1$-dimensional
        Lebesgue measure.

        \item 
        For every $t$ in the domain of $\eta$, there exist a neighborhood $U^{i,k}_{\eta,t}\ni\eta(t)$ and a smooth solution $u^{i,k}_{\eta,t}$ of the Hamilton{--}Jacobi equation $H^{i,k}(\de u^{i,k}_{\eta,t})=0$ on $U^{i,k}_{\eta,t}$, where $H^{i,k}$ is the Hamiltonian associated with $L^{i,k}$, and 
        \[
          \dot\eta=(\cL^{i,k}\,\de u^{i,k}_{\eta,t})\circ\eta
          \qquad \text{and} \qquad 
          (\div_{\mm}(\cL^{i,k}\,\de u^{i,k}_{\eta,t}))\circ\eta=-\dot\psi^{i,k}_\eta
        \]
        on $\eta^{-1}(U^{i,k}_{\eta,t})$, where $\cL^{i,k}$ is the Legendre transform associated with $L^{i,k}$.\footnote{In~\cite{Assouline24}, the function $\psi^{i,k}_\eta$ and the solution $u^{i,k}_{\eta,t}$ are stated to be of class $C^2$ and $C^3$, respectively; since the Lagrangian $L^{i,k}$ is smooth, the same proof shows that they can be taken to be smooth.}
      \end{itemize}

      By Proposition~\ref{prop:localLagrangian}\ref{item:xkL} and Lemma~\ref{lem:strain}\ref{item:dotgamma},
       for every $y \in \cT_k\cap U^{i,k}$ the Lagrangians $L^{i,k}$ and $L$ coincide in a neighborhood of
      $\dot\gamma_y(t_y) = \cL(\de u\vert_y)$,
      whence the Hamiltonians $H^{i,k}$ and $H$, as well as the corresponding Legendre transforms, coincide on a neighborhood of $\de u\left(\cT_k\cap U^{i,k}\right)$.
      Moreover, for $\nu^{i,k}$-a.e.\ $(i,k)$-ray $\eta$ and every $t \in \eta^{-1}(\cT_k)$, the covectors $\de u^{i,k}_{\eta,t}$ and $\de u$ coincide at the point $y:=\eta(t)$, since the Legendre transform is injective and maps both of them to the vector $\dot\eta(t)=\dot\gamma_y(t_y)$ (the $(i,k)$-ray $\eta$ is a sub-arc of $\gamma_y$ by Proposition~\ref{prop:localLagrangian}\ref{item:xk-calibrated}).
      Hence, for $\nu^{i,k}$-a.e.\ $(i,k)$-ray $\eta$ and every $t \in \eta^{-1}(\cT_k)$, possibly after shrinking $U^{i,k}_{\eta,t}$, we have
      \begin{equation}
        \label{eq:hj-local}
        H(\de u^{i,k}_{\eta,t})
        = 0,
        \qquad\text{ on }
        U^{i,k}_{\eta,t}
      \end{equation}
      and
      \begin{equation}\label{eq:hj-local-ray}
          \dot\eta = (\cL\,\de u^{i,k}_{\eta,t})\circ\eta
          \quad\text{and}
          \quad
        \big(
        \div_{\mm}(\cL\,\de u^{i,k}_{\eta,t})
        \big)
        \circ\eta
        =
        -\dot\psi_\eta^{i,k}
        ,
        \quad\text{ on }
        \eta^{-1}(U^{i,k}_{\eta,t}).
      \end{equation}

      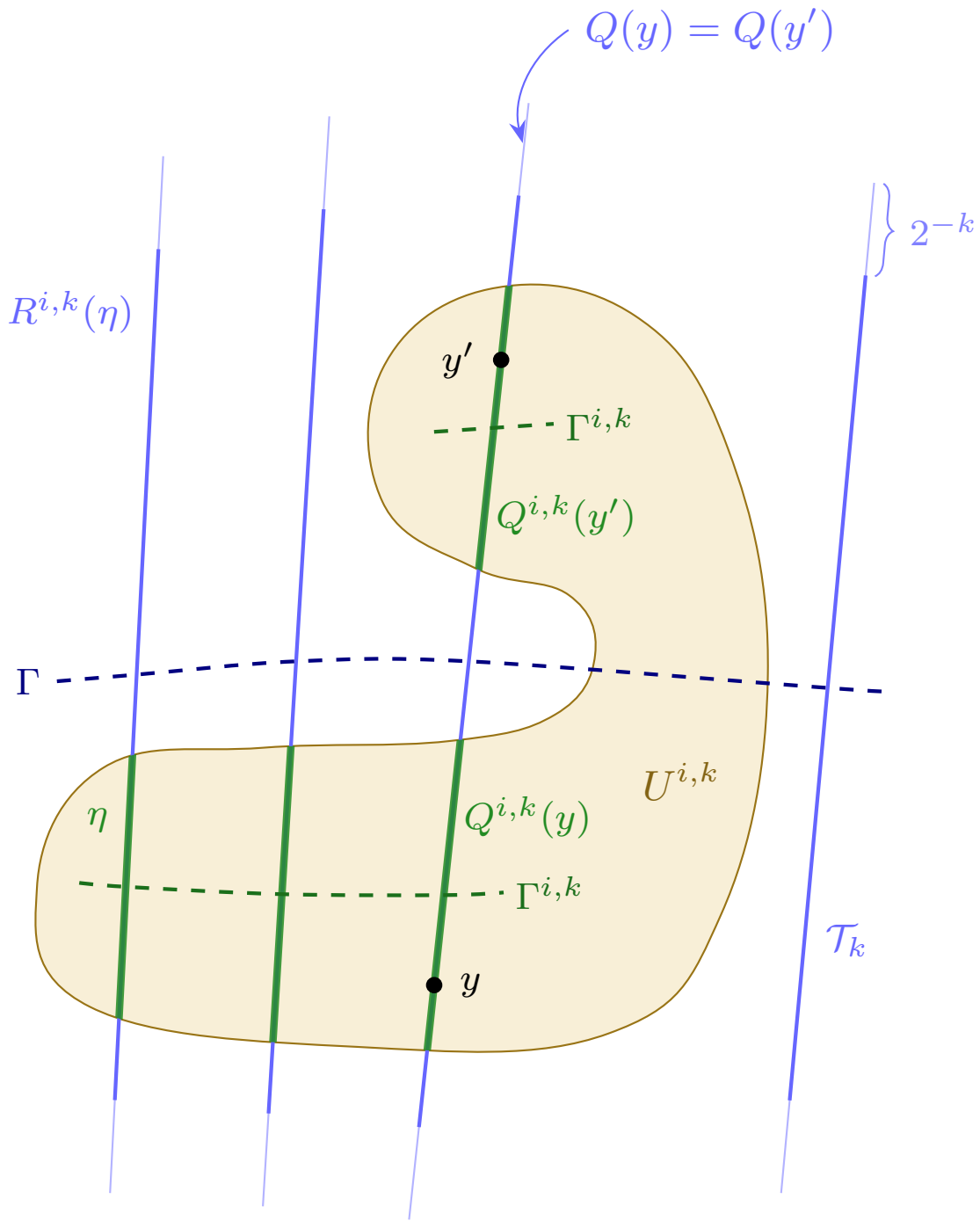
\begin{figure}[t]
      \centering
      \begin{tikzpicture}[scale=1, >=Stealth, font=\small,
          dot/.style={circle, fill, inner sep=0pt, minimum size=3.5pt}]

        \def\Upath{plot[smooth cycle, tension=.7] coordinates {
            (-2.7,2.0) (-2.3,1.1) (-0.3,0.8) (1.6,0.9) (2.45,1.8)
            (2.8,3.6) (2.5,5.3) (1.8,6.3) (0.7,6.5) (-0.1,5.9)
            (-0.1,4.9) (0.6,4.4) (1.3,4.2) (1.5,3.75) (1.2,3.3)
            (0.4,3.1) (-1.0,3.05) (-2.2,2.9)}}
        \fill[Goldenrod, fill opacity=.18] \Upath;
        \draw[Goldenrod!70!black, thin] \Upath;
        \node[Goldenrod!60!black] at (2.15,2.8) {$U^{i,k}$};

        \coordinate (P1) at (0.1,-0.5); \coordinate (Q1) at (1.0,7.9);   
        \coordinate (P2) at (-1.0,-0.4); \coordinate (Q2) at (-0.5,7.8);
        \coordinate (P3) at (2.9,-0.3); \coordinate (Q3) at (3.6,7.3);   
        \coordinate (P4) at (-2.15,-0.3); \coordinate (Q4) at (-1.75,7.5);

        \foreach \i in {1,2,3,4} {
          \draw[blue!30, thin] (P\i) -- (Q\i);
          \draw[blue!60, thick]  ($(P\i)!0.7cm!(Q\i)$) -- ($(Q\i)!0.7cm!(P\i)$);
        }
        \begin{scope}
          \clip \Upath;
          \foreach \i in {1,2,4} {
            \draw[ForestGreen, line width=1.6pt, opacity=.85] (P\i) -- (Q\i);
          }
        \end{scope}
        \node[blue!60, anchor=west] (gammalab) at ($(Q1)+(0.3,0.55)$) {$Q(y)=Q(y')$};
        \draw[->, blue!60, thin, shorten >=1pt] (gammalab.west) to[bend right=35] ($(P1)!0.96!(Q1)$);
        \node[ForestGreen, font=\footnotesize, right=-1pt] at ($(P1)!0.36!(Q1)$) {$Q^{i,k}(y)$};
        \node[ForestGreen, font=\footnotesize, right=-1pt] at ($(P1)!0.63!(Q1)$) {$Q^{i,k}(y')$};
        \node[ForestGreen, font=\footnotesize, left=1pt] at ($(P4)!0.36!(Q4)$) {$\eta$};
        \node[blue!60, font=\footnotesize, left=1pt] at ($(P4)!0.85!(Q4)$) {$\rayproj(\eta)$};
        \draw[blue!50, decorate, decoration={brace, amplitude=3pt, raise=2pt}]
          (Q3) -- ($(Q3)!0.7cm!(P3)$)
          node[midway, right=5pt, font=\footnotesize] {$2^{-k}$};
        \node[blue!60, font=\footnotesize, right=1pt] at ($(P3)!0.25!(Q3)$) {$\cT_k$};

        \path[name path=Ubdry] \Upath;
        \foreach \i in {1,2,4} { \path[name path global=ray\i] (P\i) -- (Q\i); }
        \path[name intersections={of=Ubdry and ray1, sort by=ray1, name=a}];
        \path[name intersections={of=Ubdry and ray2, sort by=ray2, name=b}];
        \path[name intersections={of=Ubdry and ray4, sort by=ray4, name=c}];
        \foreach \i in {1,2,3,4} { \coordinate (m\i) at ($(P\i)!.5!(Q\i)$); }
        \coordinate (g1) at ($(a-1)!.5!(a-2)$);
        \coordinate (h1) at ($(a-3)!.5!(a-4)$);
        \coordinate (g2) at ($(b-1)!.5!(b-2)$);
        \coordinate (g4) at ($(c-1)!.5!(c-2)$);
        \coordinate (mL) at ($(m4)+(-.6,-.05)$);
        \coordinate (mR) at ($(m3)+(.5,-.03)$);
        \coordinate (gL) at ($(g4)+(-.35,.03)$);
        \coordinate (gR) at ($(g1)+(.45,.02)$);
        \draw[NavyBlue, dashed, thick] plot[smooth, tension=.6] coordinates {(mL) (m4) (m2) (m1) (m3) (mR)};
        \node[NavyBlue, font=\footnotesize, left] at (mL) {$\Gamma$};
        \draw[ForestGreen!80!black, dashed, thick] plot[smooth, tension=.6] coordinates {(gL) (g4) (g2) (g1) (gR)};
        \draw[ForestGreen!80!black, dashed, thick] ($(h1)+(-.45,-.03)$) -- ($(h1)+(.45,.03)$);
        \node[ForestGreen!80!black, font=\footnotesize, right=-1pt] at ($(h1)+(.45,.03)$) {$\Gamma^{i,k}$};
        \node[ForestGreen!80!black, font=\footnotesize, right=-1pt] at (gR) {$\Gamma^{i,k}$};

        \coordinate (y) at ($(P1)!0.21!(Q1)$);
        \node[dot] at (y) {};
        \node[right=2pt, font=\footnotesize] at (y) {$y$};
        \coordinate (y') at ($(P1)!0.77!(Q1)$);
        \node[dot] at (y') {};
        \node[left=2pt, font=\footnotesize] at (y') {$y'$};
      \end{tikzpicture}
          \caption{The points $y,y'\in\cT_k\cap U^{i,k}$ lie on different $(i,k)$-rays $Q^{i,k}(y)\neq Q^{i,k}(y')$, but share a single transport ray $Q(y)=Q(y')$. The $(i,k)$-ray $\eta \in \Gamma^{i,k}$ is an arc of the transport ray $R^{i,k}(\eta) \in \Gamma$.}\label{fig:local-rays}
      \end{figure}

      \medskip
      We now show that the global disintegration~\eqref{eq:disintglobal}
      essentially coincides with the local one~\eqref{eq:disint-local} on the set $\cT_k\cap U^{i,k}$.
      By Proposition~\ref{prop:localLagrangian}\ref{item:xk-calibrated}, every $y\in\cT_k\cap U^{i,k}$ lies on an $(i,k)$-ray, which is a sub-arc of the transport ray $\gamma_y$. 
      Hence
      \begin{equation}\label{eq:Tk-in-Tik}
        \cT_k\cap U^{i,k}\subseteq\cT^{i,k},
      \end{equation}
      and, denoting by $\Gamma^{i,k}$ the set of $(i,k)$-rays which meet the set $\cT_k$, there is a map 
      \[
        \rayproj:\Gamma^{i,k}\to\Gamma
      \]
      sending an $(i,k)$-ray to the unique transport ray containing it.
      Writing
      \[
        Q^{i,k}:\cT_k\cap U^{i,k}\to\Gamma^{i,k}
      \]
      for the map sending each point to the $(i,k)$-ray through it,
      we thus have
      \begin{equation}\label{eq:ray-factor}
          Q = \rayproj\circ Q^{i,k}
          \qquad
          \text{on } \cT_k\cap U^{i,k},
      \end{equation}
      see Figure~\ref{fig:local-rays}.
      
      \medskip
      We endow each $(i,k)$-ray $\eta\in\Gamma^{i,k}$ with the parametrization
      it inherits as a restriction of the transport ray 
      $\gamma = \rayproj(\eta)$ to an interval $I_\eta\subseteq (0,\ell_\gamma)$.

      \medskip
      Replacing the measure $\nu^{i,k}$ by its restriction to the set $\Gamma^{i,k}$, and each measure $\mm^{i,k}_\eta$ by its restriction to the set $\cT_k$, we obtain a disintegration
      \begin{equation}\label{eq:disint-local}
        (\nu^{i,k},(\mm^{i,k}_\eta\llcorner_{\cT_k})_{\eta\in\Gamma^{i,k}})
      \end{equation}
      of the measure
      \[
        \mm^{i,k}:=\mm\llcorner_{\cT_k\cap U^{i,k}}
      \]
      with respect to $Q^{i,k}$.
      Indeed, $\cT_k \cap U^{i,k} \subseteq \cT^{i,k}$ by~\eqref{eq:Tk-in-Tik}, and by~\eqref{eq:local-needle}, for $\nu^{i,k}$-a.e.\ $(i,k)$-ray $\eta$, the measure $\mm^{i,k}_\eta\llcorner_{\cT_k}$ is nonzero if and only if $\eta\in\Gamma^{i,k}$.

      \medskip
      As we did for transport rays, we identify each $(i,k)$-ray $\eta\in\Gamma^{i,k}$ with the midpoint of the arc $\eta\cap\cT_k$, so that $\Gamma^{i,k}$ becomes a Borel subset of $M$ and the map $Q^{i,k}$ becomes Borel (for measurability see~\cite[Lem.~4.9]{Assouline24} and Lemma~\ref{lem:strain}\ref{item:borel}).
      Thus we may regard $\nu^{i,k}$ as a $\sigma$-finite Borel measure on $\Gamma^{i,k}$ satisfying
      \begin{equation}\label{eq:nu-local-ac}
        \nu^{i,k}\ll Q^{i,k}_\pushf\mm^{i,k}.
      \end{equation}
      
      Indeed, if $B\subseteq\Gamma^{i,k}$ is a Borel set, then the function $\eta\mapsto\mm^{i,k}_\eta\big((Q^{i,k})^{-1}(B)\big)$ is $\nu^{i,k}$-measurable, and the set where it is nonzero coincides with $B$ up to a $\nu^{i,k}$-null set because the measure $\mm^{i,k}_\eta\llcorner_{\cT_k}$ is nonzero and concentrated on the preimage of $\eta$ for $\nu^{i,k}$-a.e.\ $\eta\in\Gamma^{i,k}$.
      Hence $B$ is $\nu^{i,k}$-measurable, and by the disintegration formula~\eqref{eq:def-disintegration-formula}, if $Q^{i,k}_\pushf\mm^{i,k}(B)=0$ then $\nu^{i,k}(B)=0$; moreover, since images of Borel sets under $Q^{i,k}$ are analytic, every $\nu^{i,k}$-measurable set differs from a Borel set by a $\nu^{i,k}$-null set, so restricting $\nu^{i,k}$ to the Borel $\sigma$-algebra preserves its $\sigma$-finiteness and its completion.

      \medskip
      Recall that $\nu$ is a probability measure on $\Gamma$ which is mutually absolutely continuous with the measure $Q_\pushf(\mm\llcorner_\cT)$, and in particular $Q_\sharp\mm^{i,k}\ll Q_\sharp(\mm\llcorner_{\cT}) \ll \nu$.
      Hence it follows from~\eqref{eq:nu-local-ac} and~\eqref{eq:ray-factor} that $(\rayproj)_\pushf\nu^{i,k}\ll\nu$,
      so we may apply Proposition~\ref{prop:existence-disintegration}
      to the measure $\nu^{i,k}$ and to the map $\rayproj:\Gamma^{i,k}\to\Gamma$
      and obtain a disintegration 
      \begin{equation}\label{eq:disint-rayproj}
        (\nu,(\nu^{i,k}_\gamma)_{\gamma\in\Gamma})
      \end{equation}
      of $\nu^{i,k}$ with respect to $\rayproj$.
      By~\eqref{eq:ray-factor} and
      Lemma~\ref{lem:composite-disintegration}, the pair 
      $
        (\nu,(\int_{\Gamma^{i,k}}\mm^{i,k}_{\eta}\llcorner_{\cT_k}\,\nu^{i,k}_\gamma(\de\eta))_{\gamma\in\Gamma})
      $
      is a disintegration of $\mm^{i,k}$ with respect to $Q$.
      On the other hand, the pair
      \begin{equation}\label{eq:disintglobal-restricted}
        (\nu,(\mm_\gamma\llcorner_{\cT_k\cap U^{i,k}})_{\gamma\in\Gamma}),
      \end{equation}
      obtained by restricting the disintegration~\eqref{eq:disintglobal} to the set $\cT_k\cap U^{i,k}$, is another disintegration of $\mm^{i,k}$ with respect to $Q$.
      Proposition~\ref{prop:disint-uniqueness} then implies that for $\nu$-a.e.\ $\gamma\in\Gamma$ and every Borel set $A\subseteq M$,
      \begin{equation}
        \label{eq:glued-ik}
        \mm_\gamma(A\cap\cT_k\cap U^{i,k})
        =
        \int_{\Gamma^{i,k}}\mm^{i,k}_{\eta}(A\cap\cT_k)\,\nu^{i,k}_\gamma(\de\eta).
      \end{equation}

      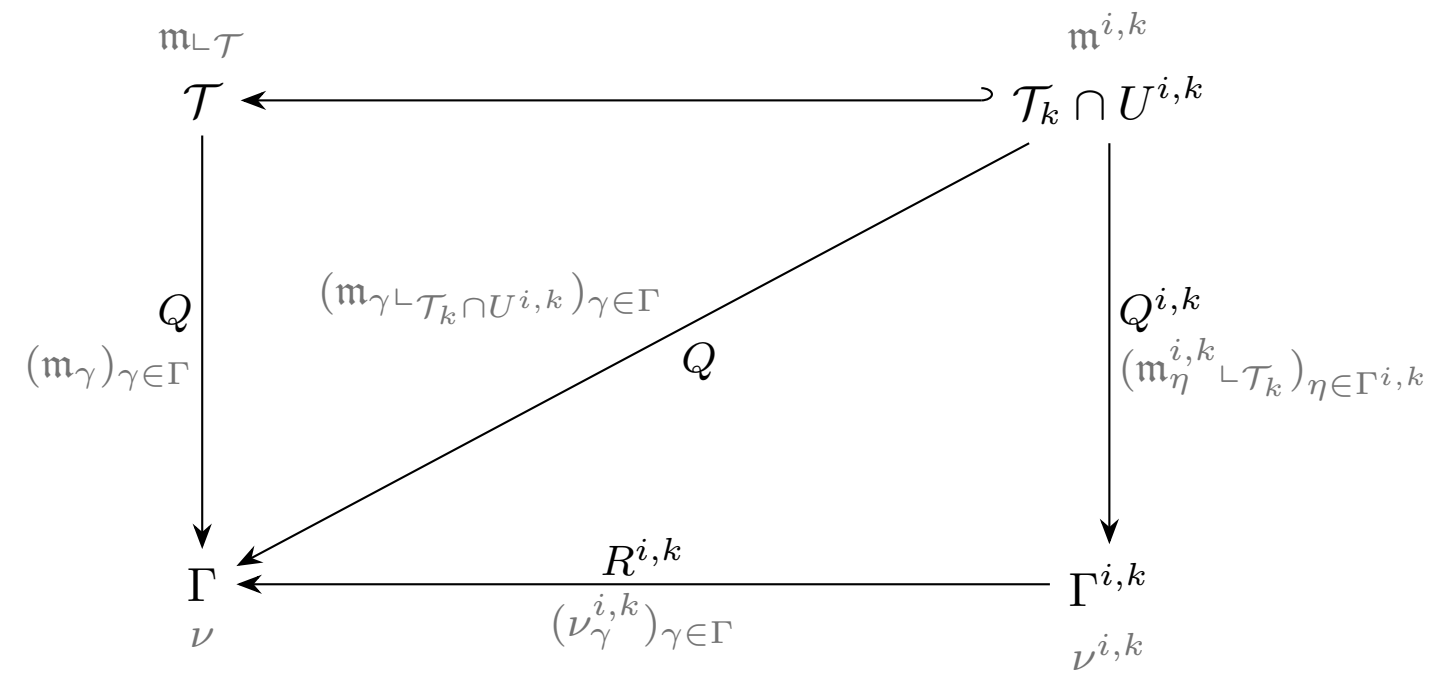
\begin{figure}[t]
      \centering
        \begin{tikzpicture}[>=Stealth, font=\small,
            inc/.style={{Hooks[left, length=2.2pt, width=5pt]}->},
            lab/.style={font=\footnotesize, inner sep=1.5pt},
            meas/.style={font=\footnotesize, text=black!55, inner sep=1.5pt}]
          \node (C) at (0,3.2)   [label={[meas]above:$\mm\llcorner_{\cT}$}] {$\cT$};
          \node (A) at (6,3.2)   [label={[meas]above:$\mm^{i,k}$}] {$\cT_k\cap U^{i,k}$};
          \node (F) at (0,0)     [label={[meas]below:$\nu$}] {$\Gamma$};
          \node (E) at (6,0)     [label={[meas]below:$\nu^{i,k}$}] {$\Gamma^{i,k}$};

          \draw[inc] (A) -- (C);
          \draw[->] (C) -- node[lab, left, align=right]
                            {$Q$\\ \textcolor{black!55}{$(\mm_\gamma)_{\gamma\in\Gamma}$}} (F);
          \draw[->] (A) -- node[lab, right, align=left]
                            {$Q^{i,k}$\\ \textcolor{black!55}{$(\mm^{i,k}_\eta\llcorner_{\cT_k})_{\eta\in\Gamma^{i,k}}$}} (E);
          \draw[->] (E) -- node[lab, above] {$\rayproj$}
                         node[lab, below] {\textcolor{black!55}{$(\nu^{i,k}_\gamma)_{\gamma\in\Gamma}$}} (F);
          \draw[->] (A) -- node[lab, above left, pos=.45]
                            {\textcolor{black!55}{$(\mm_\gamma\llcorner_{\cT_k\cap U^{i,k}})_{\gamma\in\Gamma}$}}
                         node[lab, below right, pos=.45] {$Q$} (F);
        \end{tikzpicture}
      \caption{The maps, measures, and disintegrations in the proof of Theorem~\ref{thm:main}.}\label{fig:disintegrations}
      \end{figure}

      Since $i$ and $k$ were arbitrary, there exists a set $\Gamma'\subseteq\Gamma$ of full $\nu$-measure such that for every $\gamma\in\Gamma'$, the measure $\mm_\gamma$ is nonzero and concentrated on $\gamma$, and for all $i,k\ge1$,~\eqref{eq:glued-ik} holds and~\eqref{eq:strong-consistency},~\eqref{eq:local-needle},~\eqref{eq:hj-local} and~\eqref{eq:hj-local-ray} hold for $\nu^{i,k}_\gamma$-a.e.\ $\eta\in\Gamma^{i,k}$
      (by~\eqref{eq:def-disintegration-formula}, every $\nu^{i,k}$-null subset of $\Gamma^{i,k}$ is $\nu^{i,k}_\gamma$-null for $\nu$-a.e.\ $\gamma\in\Gamma$).
      
      \medskip
      To prove~\ref{item:needle}, fix $\gamma\in\Gamma'$.
      For all $i,k\ge1$ and every $\eta\in(\rayproj)^{-1}(\gamma)$, applying~\eqref{eq:glued-ik} to subsets of the curve $\eta$ gives
      \begin{equation}\label{eq:glued-on-local-ray}
        \mm_\gamma\llcorner_{\cT_k\cap\eta}
        =
        \nu^{i,k}_\gamma\big(\{\eta\}\big)\,
        \mm^{i,k}_{\eta}\llcorner_{\cT_k},
      \end{equation}
      since for $\nu^{i,k}_\gamma$-a.e.\ $\eta'\neq\eta$, the measure $\mm^{i,k}_{\eta'}\llcorner_{\cT_k}$ is concentrated on $\eta'$ by~\eqref{eq:strong-consistency}, and $\eta'$ is disjoint from $\eta$.
      For every $t\in(0,\ell_\gamma)$ there exist $i,k\ge1$ such that $\gamma(t)\in\cT_k\cap U^{i,k}$, and by~\eqref{eq:Tk-in-Tik} and~\eqref{eq:ray-factor}, the $(i,k)$-ray $\eta:=Q^{i,k}(\gamma(t))$ satisfies $\rayproj(\eta)=\gamma$ and $t\in I_\eta$.
      Thus the open intervals 
      \[
        \gamma^{-1}(\cT_k)\cap I_\eta=I_\eta\cap(2^{-k},\ell_\gamma-2^{-k}), \qquad i,k\ge1, \quad \eta\in(\rayproj)^{-1}(\gamma),
      \]
      cover the interval $(0,\ell_\gamma)$.
      By~\eqref{eq:glued-on-local-ray} and~\eqref{eq:local-needle}, on each of these intervals the measure $m_\gamma:=(\gamma^{-1})_\pushf \mm_\gamma$ has a density with respect to $\Leb^1$ which is either identically zero, or smooth, positive and proportional to $e^{-\psi^{i,k}_{\eta}}$.
      But these open intervals cover the connected interval $(0,\ell_\gamma)$ and $m_\gamma\neq0$, so the density is smooth and positive on $(0,\ell_\gamma)$,
      i.e., we can write $m_\gamma=e^{-\psi_\gamma}\Leb^1$
      where $\psi_\gamma:(0,\ell_\gamma)\to \RR$ is a smooth function. Since this density is positive,~\eqref{eq:glued-on-local-ray} gives $\nu^{i,k}_\gamma(\{\eta\})>0$ for all $i,k\ge1$ and every $\eta\in(\rayproj)^{-1}(\gamma)$; hence~\eqref{eq:local-needle} holds for every such $\eta$ by the choice of $\Gamma'$, and
      $\psi_\gamma-\psi^{i,k}_\eta$ is constant on $\gamma^{-1}(\cT_k)\cap I_\eta$ for all $i,k\ge1$ and every $\eta\in(\rayproj)^{-1}(\gamma)$.
      Since the measure $\mm_\gamma$ is concentrated on $\gamma$, this proves~\ref{item:needle} for all $\gamma\in\Gamma'$.
 
      \medskip
      Assertion~\ref{item:massbalance} follows
      from Lemma~\ref{lem:strain}\ref{item:balance} by a standard
      argument, see the proofs of~\cite[Th.~1.5(B)]{Kl17}
      and~\cite[Th.~5.9]{Assouline24}.
      
      \medskip
      It remains to prove~\ref{item:HJneedle}.
      Let $\gamma\in\Gamma'$ and $t\in(0,\ell_\gamma)$, and choose $i,k\ge1$ and $\eta\in(\rayproj)^{-1}(\gamma)$ such that $t\in\gamma^{-1}(\cT_k)\cap I_\eta$.
      By the choice of $\Gamma'$,~\eqref{eq:hj-local} and~\eqref{eq:hj-local-ray} hold $\nu^{i,k}_\gamma$-almost everywhere on $\Gamma^{i,k}$, hence they hold at $\eta$ since $\nu^{i,k}_\gamma(\{\eta\})>0$.
      Set $u_{\gamma,t}:=u^{i,k}_{\eta,t}$, so that $H(\de u_{\gamma,t})=0$ on $U^{i,k}_{\eta,t}$ by~\eqref{eq:hj-local}; since $\psi_\gamma-\psi^{i,k}_\eta$ is constant on $\gamma^{-1}(\cT_k)\cap I_\eta$,~\eqref{eq:HJneedle} holds on the open set $\gamma^{-1}(U^{i,k}_{\eta,t})\cap\gamma^{-1}(\cT_k)\cap I_\eta \ni t$.
      \qed

    \section{Displacement interpolations}\label{sec:di}
      In this section we use Theorem~\ref{thm:main} to construct 
      displacement interpolations, which are dynamical solutions to the optimal transport problem
      for the cost $\cc$.
      Recall that
      \[
            SM : = E^{-1}(0)\subseteq\interior\dom L.
      \]

      Denote by $\Phi$ the restriction of the Euler-Lagrange flow
      to $SM$, so that 
      \[
        \Phi(v,\tau) = \dot\gamma_v(\tau),
        \qquad
        v \in SM,
      \]
      where $\gamma_v$ is
      the solution to~\eqref{eq:EL} with initial velocity
      $v$. Denote by $\mmz$ the closed subset of
      $SM\times{[}0,\infty)$ consisting of pairs $(v,\ell)$ for which
      the curve $\gamma_v(\tau) = \pi\left(\Phi(v,\tau)\right)$ is
      action-minimizing on ${[}0,\ell]$, where $\pi:SM \to M$ is the natural projection, i.e.
      \[
        \mmz : =
        \left\{
          (v,\ell) \in SM \times {[}0,\infty)
          \,\, \colon \, \,
          \cc\left(\gamma_v(0),\gamma_v(\ell)\right)
          =
          \action\left[\gamma_v\vert_{[0,\ell]}\right]
        \right\}.
      \]
      
      Note that elements of $\mmz$ are in one-to-one correspondence with minimizing extremals defined on intervals of the form $[0,\ell]$ for $\ell \ge 0$. 

      \begin{definition}[\text{Initial transport measures, cf.~\cite[Sect.~2]{BernardBuffoni}}]\label{def:initial}
          Let $\cX_0\subseteq M$ be a Borel set and let $\mu_0$ be a Borel measure on $\cX_0$.
          An \emph{initial transport measure} for $\mu_0$ is a measure on $\mmz$ of the form 
          \[
              \tilde \mu_0  : = (V_0,\ell)_\sharp\mu_0
          \]
          where $V_0:\cX_0 \to SM$ is a Borel vector field, i.e., a section of $SM$:
          \[
              \pi\circ V_0 = \mathrm{id},
          \]
          and $\ell:\cX_0 \to [0,\infty)$ is a $\mu_0$-measurable function.
          The initial transport measure $\tilde\mu_0$ is called \emph{monotone} if for all $x_0,x_0' \in \cX_0$ and all $\tau \ge 0$,
          \[
              \Phi(V_0\vert_{x_0},\tau) = V_0\vert_{x_0'}
              \qquad \implies \qquad
              \ell(x_0) \le \ell(x_0') + \tau.
          \]

          Monotonicity means that the initial transport measure respects the order of points along each trajectory of the flow $\Phi$: if the Euler{--}Lagrange trajectory issuing from $x_0$ at velocity $V_0\vert_{x_0}$ passes through $x_0'$ after time $\tau$, then it must terminate no later than the trajectory issuing from $x_0'$ at velocity $V_0\vert_{x_0'}$.
      \end{definition}

      \begin{definition}[Displacement interpolations]\label{def:di}
          A family $(\mu_\lambda)_{0 \le \lambda \le 1}$ of Borel probability measures on $M$ is called a \emph{displacement interpolation} if there exists an initial transport measure $\tilde\mu_0 = (V_0,\ell)_\sharp\mu_0$ for $\mu_0$ such that for all $\lambda \in [0,1]$
          \begin{equation}\label{eq:Tlambda-def}
              \mu_\lambda = (\rT_\lambda)_\sharp\mu_0
              \qquad\text{where}\qquad
              \rT_\lambda : = \pi\circ\Phi\circ(V_0,\lambda\ell) : \cX_0 \to M,
          \end{equation}
          and such that the map
          \[
              \rT : = \rT_1 = \pi\circ\Phi\circ(V_0,\ell) : \cX_0 \to M
          \]
          is an optimal transport map between the endpoint measures $\mu_0$ and $\mu_1 = \rT_\sharp\mu_0$.
          We then define a family of measures $(\dot\mu_\lambda)_{0 \le \lambda \le 1}$ on $SM$ by
          \[
              \dot\mu_\lambda : = (\Phi\circ(V_0,\lambda\ell))_\sharp\mu_0, \qquad \lambda \in [0,1].
          \]

          Note that
          \[
              \pi_\sharp\dot\mu_\lambda = \mu_\lambda, \qquad \lambda \in [0,1].
          \]

          We refer to the function $\ell$ as the \emph{duration function} of the displacement interpolation $(\mu_\lambda)$. A displacement interpolation is called \emph{monotone} if the initial transport measure $\tilde\mu_0$ may be chosen monotone.
      \end{definition}

        \begin{theorem}[Existence and uniqueness of monotone displacement interpolations]\label{thm:diexistence}
          Let $M$ be a smooth manifold,
          let $L:TM\to\RR\cup\{+\infty\}$ be a Lagrangian satisfying assumptions~\ref{item:convex}{--}\ref{item:positivezerosection} of Section~\ref{sec:Lagrangians},
          and let $\cA \subseteq (M\times M)\setminus\Delta$
          satisfy assumptions~\ref{item:attained} and~\ref{item:compactness}.
          Let $\mm$ be a measure on $M$ with a smooth positive density,
          and let $\mu_0=f_0\mm$ and $\mu_1=f_1\mm$ be compactly supported,
          absolutely continuous Borel probability measures on $M$ such that $\supp\mu_0\times\supp\mu_1 \subseteq \cA$.

          \medskip
          Let $\nu$ and $(m_\gamma)_{\gamma\in\Gamma}$ be as in the conclusion of Theorem~\ref{thm:main}, and write
          \[
            m_{\gamma,i} := (f_i\circ\gamma)\,m_\gamma, \qquad i = 0,1.
          \]

          For $\nu$-a.e.\ $\gamma \in \Gamma$,
          let $(m_{\gamma,\inter})_{0 \le \inter \le 1}$ be the monotone displacement interpolation between $m_{\gamma,0}$ and $m_{\gamma,1}$ on the real line (see Remark~\ref{rem:1d-di-exists}).
          
          \medskip
          Then the family
          \begin{equation}\label{eq:mulambdadisint}
            \mu_\inter
            :=
            \int_\Gamma \gamma_{\pushf}(m_{\gamma,\inter})\,\nu(\de\gamma),
            \qquad \inter \in [0,1],
          \end{equation}
          is a monotone displacement interpolation between $\mu_0$ and $\mu_1$, and the measures $(\dot\mu_\inter)_{0\le\inter\le1}$ in Definition~\ref{def:di} are given by
          \begin{equation}\label{eq:dotmudisint}
            \dot\mu_\inter
            =
            \int_\Gamma \dot\gamma_{\pushf}(m_{\gamma,\inter})\,\nu(\de\gamma),
            \qquad \inter \in [0,1],
          \end{equation}
          where $\dot\gamma : (0,\ell_\gamma) \to SM$ denotes the velocity of $\gamma$.
          Each $\mu_\inter$ is absolutely continuous with respect to $\mm$, and its density $f_\inter$ satisfies
          \begin{equation}\label{eq:givenDI-densities}
            m_{\gamma,\inter} = (f_\inter\circ\gamma)\,m_\gamma
            \qquad\text{for $\nu$-a.e.\ $\gamma\in\Gamma$.}
          \end{equation}

          \medskip
          Moreover, $(\mu_\inter)_{0\le\inter\le1}$ is the unique monotone displacement interpolation between $\mu_0$ and $\mu_1$, and any two pairs $(V_0,\ell)$ as in Definition~\ref{def:di} whose induced initial transport measures are monotone agree $\mu_0$-almost everywhere. In particular, \emph{every} monotone displacement interpolation between $\mu_0$ and $\mu_1$ is given by~\eqref{eq:mulambdadisint}.
        \end{theorem}
        \begin{remark}\label{rem:1d-di-exists}
          By~\eqref{eq:mass-balance-main-theorem}, for $\nu$-a.e.\ $\gamma\in\Gamma$ the measures $m_{\gamma,0}$ and $m_{\gamma,1}$ have equal finite mass, hence the monotone displacement interpolation between them exists.
        \end{remark}
        \begin{proof}
          We freely use the notation of Section~\ref{sec:OT}.
          By Theorem~\ref{thm:main}\ref{item:disint}--\ref{item:massbalance}, for $\nu$-a.e. $\gamma \in \Gamma$,
          the measures $m_{\gamma,0}$ and $m_{\gamma,1}$ on the interval $(0,\ell_\gamma)$ are absolutely continuous, of equal finite mass, and satisfy
          \begin{equation}\label{eq:tail-condition}
            m_{\gamma,0}\left((t,\ell_\gamma)\right) \le m_{\gamma,1}\left((t,\ell_\gamma)\right) \qquad \text{for all $t \in (0,\ell_\gamma)$.}
          \end{equation}
          Since $\mu_0\perp\mu_1$, we also have $m_{\gamma,0}\perp m_{\gamma,1}$ for $\nu$-a.e.\ $\gamma\in\Gamma$ by~\eqref{eq:disint}.
          
          By~\eqref{eq:disint}, removing a $\nu$-negligible set from $\Gamma$ does not change the fact that $\mu_0(\cT) = \mu_1(\cT) = 1$.
          Thus we may assume that the above properties hold for every $\gamma \in \Gamma$.

          \medskip
          For $\gamma \in \Gamma$, let $\rT_\gamma : (0,\ell_\gamma) \to (0,\ell_\gamma]$ be the nondecreasing map satisfying $(\rT_\gamma)_\pushf\, m_{\gamma,0} = m_{\gamma,1}$,
          see~\cite[Ch.~2]{Santambrogio15} or~\cite[Lem.~5.12]{Assouline24}.
          By~\eqref{eq:tail-condition} we have $\rT_\gamma(t) \ge t$ for $m_{\gamma,0}$-a.e. $t \in (0,\ell_\gamma)$; replacing $\rT_\gamma$ by $\max\{\rT_\gamma,\id\}$, which is still nondecreasing, we may assume that $\rT_\gamma \ge \id$ everywhere.
          Set $\rT_{\gamma,\inter} := (1-\inter)\id + \inter\,\rT_\gamma$ for $\inter \in [0,1]$, so that
          \begin{equation}\label{eq:mgammalambdadef}
            m_{\gamma,\inter} = (\rT_{\gamma,\inter})_\pushf\, m_{\gamma,0}.
          \end{equation}
          
          Define
          \begin{equation}\label{eq:V0elldef}
              V_0\vert_x := \dot\gamma_x(t_x),
              \qquad
              \ell(x) := \rT_{\gamma_x}(t_x) - t_x \ge 0,
              \qquad
              x \in \cT.
          \end{equation}

          Since transport rays are minimizing extremals, the vector field $V_0$ takes values in $SM$ (see Section~\ref{sec:Lagrangians}), and $V_0 = \cL\circ\de u$ on $\cT$ by Lemma~\ref{lem:strain}~\ref{item:dotgamma}; in particular, $V_0$ is Borel by Lemma~\ref{lem:strain}~\ref{item:lipschitz}.
          The verification that $\ell$ is $\mu_0$-measurable is routine, cf.~\cite[Lem.~2.9]{Kl17}.
          Extend $(V_0,\ell)$ to $\cX_0$ by setting $\ell := 0$ on $\cX_0\setminus\cT$ and extending $V_0$ in an arbitrary Borel fashion.
          
          \medskip
          Since $\mu_0\perp\mu_1$, we have $\ell(x) > 0$ for $\mu_0$-a.e.\ $x \in \cT$.
          For such $x$, the curve
          \[
              \tau \mapsto \pi\left(\Phi(V_0\vert_x,\tau)\right) = \gamma_x(t_x + \tau)
          \]
          is a restriction of the minimizing extremal $\bar\gamma_x$, hence action-minimizing on $[0,\ell(x)]$; that is, $(V_0\vert_x,\ell(x)) \in \mmz$. Therefore $\tilde\mu_0 := (V_0,\ell)_\sharp\mu_0$ is an initial transport measure for $\mu_0$ in the sense of Definition~\ref{def:initial}.

          \medskip
          Let $\rT_\inter$ be as in~\eqref{eq:Tlambda-def}.
          Then, if $\mu_\lambda$ is defined via~\eqref{eq:mulambdadisint}, it follows from~\eqref{eq:m-is-a-push-forward},~\eqref{eq:disint} and~\eqref{eq:mgammalambdadef} that $(\rT_\inter)_\pushf\,\mu_0 = \mu_\inter$ for all $\inter \in [0,1]$ and that $\dot\mu_\lambda$ is given by~\eqref{eq:dotmudisint}.
          For $\mu_0$-a.e.\ $x$, the points $x\in\cX_0\cap\cT$ and $\rT(x)\in\cX_1\cap\cT$ lie on the transport ray $\gamma_x$, which is calibrated by Lemma~\ref{lem:strain}\ref{item:u-facts}, so that $\phi_1(\rT(x))-\phi_0(x) = u(\rT(x))-u(x) = \cc(x,\rT(x))$ by~\eqref{eq:phi-u-coincide}.
          By Proposition~\ref{prop:KD}, applied to the coupling $(\id,\rT)_\sharp\mu_0$, the map $\rT$ is an optimal transport map between $\mu_0$ and $\mu_1$.
          Hence the family $(\mu_\inter)_{0\le\inter\le1}$ is a displacement interpolation between $\mu_0$ and $\mu_1$ in the sense of Definition~\ref{def:di}.
          Monotonicity of the initial transport measure $\tilde\mu_0$ and, thus, of the displacement interpolation $(\mu_\inter)_{0\le\inter\le1}$, follows from monotonicity of the maps $\rT_\gamma$.

          \medskip
          Each measure $m_{\gamma,\inter}$ is absolutely continuous;
          since the measures $m_\gamma$ in Theorem~\ref{thm:main}~\ref{item:needle} have positive densities, the disintegration~\eqref{eq:disint} then shows that the measure $\mu_\inter$ given by~\eqref{eq:mulambdadisint} is absolutely continuous with respect to $\mm$ for every $\inter \in [0,1]$, and~\eqref{eq:givenDI-densities} follows from~\eqref{eq:disint} and~\eqref{eq:mulambdadisint}.

          \medskip
          It remains to prove uniqueness. Let $(\tilde V_0,\tilde\ell)$ denote an arbitrary pair as in Definition~\ref{def:di}
          whose induced initial transport measure $(\tilde V_0,\tilde\ell)_\sharp\mu_0$
          is monotone according to Definition~\ref{def:initial},
          with maps $\tilde\rT_\inter = \pi\circ\Phi\circ(\tilde V_0,\inter\tilde\ell)$ and $\tilde\rT := \tilde\rT_1$,
          inducing a displacement interpolation $(\mu'_\inter)_{0\le\inter\le1}$.
          It suffices to show that $(\tilde V_0,\tilde\ell) = (V_0,\ell)$ $\mu_0$-a.e., where $(V_0,\ell)$ is the pair constructed in~\eqref{eq:V0elldef}.

          \medskip
          Since $\tilde\rT$ is an optimal transport map and $(u,u)$ is a Kantorovich pair, Proposition~\ref{prop:KD} implies that
          \begin{equation}\label{eq:Tcontact}
              u(\tilde\rT(x)) - u(x) = \cc(x,\tilde\rT(x)) \qquad \text{for $\mu_0$-a.e.\ $x \in \cX_0$.}
          \end{equation}

          Since $\mu_0\perp\mu_1$ and the set $\cT$ has full $\mu_0$-measure, we also have that
          \begin{equation}\label{eq:aeonX0}
            x \in \cX_0\cap\cT,
            \quad
            \tilde\ell(x) > 0,
            \quad
            \text{and} \quad
            (\tilde V_0\vert_x,\tilde\ell(x)) \in \mmz
            ,
            \qquad \text{for $\mu_0$-a.e. $x \in \cX_0$.}
          \end{equation}
          
          Fix $x \in \cX_0\cap\cT$ satisfying~\eqref{eq:Tcontact} and~\eqref{eq:aeonX0}.
          Let
          \[
            \sigma(s) := \pi\left(\Phi(\tilde V_0\vert_x,s)\right), \qquad s \in [0,\tilde\ell(x)].
          \]
          By~\eqref{eq:aeonX0}, $\sigma$ is a minimizing extremal joining $x$ to $\tilde\rT(x)$.
          By~\eqref{eq:Tcontact}, $u(\tilde\rT(x)) - u(x) = \cc(x,\tilde\rT(x)) = \action[\sigma]$.
          Since $u$ is $L$-dominated, it follows that $u(\sigma(s)) = u(x) + \action[\sigma\vert_{{[0,s]}}]$ for $0 \le s \le \tilde\ell(x)$.
          Differentiating at $s = 0$ ($u$ is differentiable at $x$ by Lemma~\ref{lem:strain}) yields $\de u\vert_x(\tilde V_0\vert_x) = L(\tilde V_0\vert_x)$; as $H(\de u\vert_x) = 0$, the vector $\tilde V_0\vert_x$ attains the supremum in~\eqref{eq:Hdef} with $p = \de u\vert_{x}$.
          Thus, by Lemma~\ref{lem:strain}~\ref{item:dotgamma} and~\eqref{eq:V0elldef},
          \begin{equation}\label{eq:V0unique}
              \tilde V_0\vert_x = \cL\left(\de u\vert_x\right) = \dot{\gamma}_x(t_x) = V_0\vert_x \qquad \text{for $\mu_0$-a.e.\ $x$.}
          \end{equation}

          \medskip
          Since $u$ is dominated and the transport ray $\gamma_x$ is
          calibrated,~\eqref{eq:Tcontact} and~\eqref{eq:aeonX0} imply that
          $\cc\left(\gamma_x(s),x\right) + \cc(x,\tilde\rT(x)) = \cc(\gamma_x(s),\tilde\rT(x)) = u(\tilde\rT(x)) - u(\gamma_x(s))$ for $0<s<t_x$.
          Since $\gamma_x$ is maximal, it follows that $t_x + \tilde\ell(x) \le \ell_{\gamma_x}$.
          Since $\sigma$ and the curve $\tau\mapsto\gamma_x(t_x+\tau)$ solve~\eqref{eq:EL} with the same initial velocity by~\eqref{eq:V0unique}, it follows that
          \begin{equation}\label{eq:Talongray}
              \tilde\rT(x) = \sigma(\tilde\ell(x)) = \gamma_x(t_x+\tilde\ell(x)).
          \end{equation}

          \medskip
          For $\nu$-a.e.\ $\gamma$, the map 
          \[
            \tilde\rT_\gamma(\tau) := \tau + \tilde\ell(\gamma(\tau))
          \]
          is defined $m_{\gamma,0}$-a.e.\ and nondecreasing by~\eqref{eq:V0unique} and monotonicity of $(\tilde V_0,\tilde\ell)$.
          The map $\tilde\rT$ maps each transport ray into itself by~\eqref{eq:Talongray}, and distinct transport rays are disjoint by Lemma~\ref{lem:strain}\ref{item:transportray}; hence~\eqref{eq:disint} and~\eqref{eq:m-is-a-push-forward} imply that
          \[
            (\tilde\rT_\gamma)_\pushf\, m_{\gamma,0} = m_{\gamma,1} \qquad \text{for $\nu$-a.e.\ $\gamma \in \Gamma$.}
          \]
          
          Extending $\tilde\rT_\gamma$ to a nondecreasing map on $(0,\ell_\gamma)$ and using the uniqueness of one dimensional monotone displacement interpolations (see~\cite[Lem.~5.12]{Assouline24}), we conclude that $\tilde\rT_\gamma = \rT_\gamma$ $m_{\gamma,0}$-a.e., i.e.\ 
          \begin{equation}\label{eq:ellunique}
            \tilde\ell(x) = \rT_{\gamma_x}(t_x) - t_x = \ell(x) \qquad \text{ for $\mu_0$-a.e.\ $x$.}
          \end{equation}
          
          By~\eqref{eq:V0unique} and~\eqref{eq:ellunique}, $(\tilde V_0,\tilde\ell)=(V_0,\ell)$ $\mu_0$-a.e., which proves uniqueness.
        \end{proof}

    \section{Application: electromagnetic optimal transport}\label{sec:em}
    
      In this section we apply the results of the previous sections to electromagnetic Lagrangians
      on Lorentzian manifolds.
      Let $(M,g)$ be an $n$-dimensional
      Lorentzian manifold with signature
      $(+ -\,\cdots\,-)$, which we assume to be \emph{time-oriented}, i.e., to admit a smooth global timelike vector field $Z$.
      
      \medskip
      A tangent vector $v \in TM$ is said to be \emph{timelike} (resp.\ \emph{causal}) if $\langle v,v\rangle > 0$ (resp.\ $\langle v,v\rangle\ge 0$), and \emph{future-directed} if it is causal and $\langle Z,v\rangle > 0$.
      A piecewise-$C^1$ curve $\gamma:[0,\ell] \to M$ is said to be
      \emph{chronological} (resp.\ \emph{causal})
      if $\dot\gamma(\tau)$ is timelike (resp.\ causal) and
      future-directed for all $\tau\in [0,\ell]$.
      For $x,y\in M$, write $x\ll y$ if there exists a chronological curve joining $x$ to $y$.
      We refer the reader to~\cite{ONeill83,BeemEhrlichEasley96} for
      more details on Lorentzian geometry.

      \medskip
      Let $A \in \Omega^1(M)$ be a one-form on $M$, and let $F := \de A$.
      We think of $A$ as an electromagnetic potential and of $F$ as the corresponding electromagnetic field.
      We denote by $F^\sharp$ the $(1,1)$-tensor given by
      \[
        g(X,F^\sharp(Y))=F(X,Y).
      \]
      Fix $r <1$, $r\neq 0$, and
      define a Lagrangian $L_A:TM \to \RR\cup\{+\infty\}$ by
      \begin{equation}\label{eq:Ldef}
          L_A({v}) :=
          \begin{cases}
            - \frac{1}{r}
            \langle v,v\rangle^{r/2} +
            \frac{1-r}{r} - A(v)
            ,
            \qquad& \text{ if $v$ is future-directed or $v = 0$,} \\
            \infty,
                  & \text{ otherwise.}
          \end{cases}
        \end{equation}
      For $r<0$, we also set $L_A(v):=\infty$ if $v$ is null or $v=0$.

      \medskip
      We will denote by $\action_A$ (resp.\ $\cc_A$) the
      action (resp.\ cost) induced by $L_A$.
      A future-directed, timelike $C^1$ curve $\gamma$ solves the Euler{--}Lagrange equation for $L_A$ if and only if it has constant speed and satisfies $|\dot\gamma|^{r-2}\nabla_{\dot\gamma}\dot\gamma = F^\sharp\dot\gamma$ (the proof is a straightforward modification of the derivation of the Riemannian geodesic equation~\cite{CheegerEbin75}).
      Since the zeroth energy level of $L_A$ consists of unit vectors, minimizing extremals of the Lagrangian $L_A$ are parametrized by proper time, and solve the Lorentz force
      equation
      \begin{equation}
          \label{eq:euler-lagrange-any-speed}
          \nabla_{\dot\gamma}{\dot\gamma}
          =
          F^\sharp \dot\gamma
          .
      \end{equation}
      
      Moreover, since $\xi-\xi^r/r\ge 1-1/r$ for all $\xi>0$, with equality precisely at $\xi=1$, we have $L_A(v)\ge -|v|-A(v)$ for every future-directed timelike vector $v$, with equality precisely when $|v|=1$.
      It follows that
      \begin{equation}\label{eq:cA-param-inv}
        \cc_A(x,y)
        =
        \inf_\gamma\action^1_A[\gamma],
        \qquad
        \action^1_A[\gamma] : =
        \int\bigl(-|\dot\gamma| - A(\dot\gamma)\bigr),
      \end{equation}
      where the infimum is over chronological curves joining $x$ to $y$; in particular, the cost $\cc_A$ does not depend on $r$.

      \medskip
      It is straightforward to verify that the Lagrangian $L_A$ satisfies assumptions~\ref{item:convex}--\ref{item:positivezerosection}.
      We now state sufficient conditions which 
      guarantee assumptions~\ref{item:attained} and~\ref{item:compactness} for a set $\cA$.
      Fix a Riemannian metric $h$ on $M$, and denote by $\Len_h$ the length with respect to $h$; in what follows, constants implied by big-O notation may depend on $g$, $A$, $h$ and on the compact set $\Omega$ below.
      For every compact set $\Omega \subseteq M$, every solution $\gamma:[0,\ell]\to \Omega$
      of~\eqref{eq:euler-lagrange-any-speed}
      satisfies $\big|\frac{\de}{\de t}\log|\dot\gamma|_h\big| = O\big(|\dot\gamma|_h + 1\big)$ and therefore
      \begin{equation}\label{eq:emvelocityestimate}
        |\dot\gamma(t)|_h
        \le
        |\dot\gamma(t')|_h
        \cdot
        \exp\big(O(\Len_h(\gamma) + \ell)\big)
        ,
        \qquad
        t,t' \in [0,\ell].
      \end{equation}

      Denote by $\I := \{(x,y)\in M\times M \colon x\ll y\}$ the chronological relation, and set
      \begin{equation}\label{eq:I-A-definition}
        \I_A : =
        \Big\{
        (x,y) \in \I \,\colon\,
        \inf_{\gamma_0}\Big(-\int_{\gamma_0}A\Big) > \cc_A(x,y)
        \Big\}
        ,
      \end{equation}
      where the infimum is over null geodesics $\gamma_0$ joining $x$ to $y$, with the convention $\inf\varnothing := +\infty$.
      Note that $-\int_{\gamma_0}A = \action^1_A[\gamma_0]$ for $\gamma_0$ null; thus $\I_A$ is the set of pairs $(x,y)$ for which the cost $\cc_A(x,y)$
      is \emph{not} attained by a null geodesic.

      \begin{remark}
      If no null geodesic joins a pair of points of a compact set
      $K \subseteq \I$, then $K\subseteq\I_A$ for any
      one-form $A$, or, in other words,
      one may take $s_0 = \infty$ in
      Lemma~\ref{lem:smallcmr} and hence in Theorem~\ref{thm:em}.
      This is the case when $M$ is strongly causal and $K$ is
      sufficiently close to the diagonal;
      more precisely, if $M$ is strongly causal, then there exists
      an open set $U \subseteq M\times M$ containing the diagonal such that
      $\I_A\cap U = \I\cap U$ for any one-form $A$.
      In Minkowski spacetime one can take $U = \RR^n \times \RR^n$.
      \end{remark}

      \begin{lemma}\label{lem:em-limit}
        Assume that $M$ is globally hyperbolic, and let $K \subseteq \I_A$ be a compact set.
        \begin{enumerate}
          \item\label{item:em-speed} There exists $C > 0$ such that every chronological solution $\gamma:[0,\ell]\to M$ of~\eqref{eq:euler-lagrange-any-speed} parametrized by proper time with endpoints in $K$ satisfies $|\dot\gamma|_h \le C/\ell$.
          \item\label{item:em-limit} Let $\gamma_n:[0,\ell_n]\to M$ be a sequence of
        chronological solutions of~\eqref{eq:euler-lagrange-any-speed}
        parametrized by proper time, such that $(\gamma_n(0),\gamma_n(\ell_n)) \in K$ converges to some $(x,y) \in K$ and $\limsup_n\action_{A}[\gamma_n] \le \cc_{A}(x,y)$.
          Then, after passing to a subsequence, $\ell_n$ converges to a positive limit, and $x$ and $y$ are joined by a minimizing extremal.
        \end{enumerate}
      \end{lemma}
      \begin{proof}
        In this proof, constants implied by big-O notation may depend also on $K$.
        By global hyperbolicity, every chronological solution $\gamma:[0,\ell]\to M$ of~\eqref{eq:euler-lagrange-any-speed} parametrized by proper time with endpoints in $K$ satisfies $\ell+\Len_h(\gamma)=O(1)$, and hence, by~\eqref{eq:emvelocityestimate},
        \[
          \max_{[0,\ell]}|\dot\gamma|_h
          =
          O\Big(\min_{[0,\ell]}|\dot\gamma|_h\Big)
          =
          O\Big(\frac{\Len_h(\gamma)}{\ell}\Big)
          =
          O\Big(\frac{1}{\ell}\Big).
        \]
        This proves~\ref{item:em-speed}.

        \medskip
        Now let $(\gamma_n)$ be as in~\ref{item:em-limit}.
        For each $n \ge 1$, let $\sigma_n$ denote the reparametrization of $\gamma_n$ by constant $g$-speed on the interval $[0,1]$. Then the curves $\sigma_n$ solve~\eqref{eq:euler-lagrange-any-speed} with $F$ replaced by $\ell_n F$, and their $h$-speeds are uniformly bounded by~\ref{item:em-speed}.
        Hence the velocities $\dot\sigma_n(0)$ lie in a compact subset of $TM$, so after passing to a subsequence, $\ell_n \to \ell$ and $\dot\sigma_n(0)$ converges; by continuous dependence of solutions of~\eqref{eq:euler-lagrange-any-speed} on initial data and on the two-form $F$, the curves $\sigma_n$ converge in $C^1$ to a curve $\sigma:[0,1]\to M$ joining $x$ to $y$, which solves~\eqref{eq:euler-lagrange-any-speed} with $F$ replaced by $\ell F$, has constant speed $\ell$, and satisfies
        \[
          -\ell - \int_{\sigma}A
          =
          \lim_n\Big(-\ell_n - \int_{\gamma_n}A\Big)
          =
          \lim_n \action_{A}[\gamma_n]
          \le
          \cc_{A}(x,y)
          .
        \]

        \medskip
        If $\ell = 0$, then $\sigma$ is a null geodesic joining $x$ to $y$, which is nonconstant since $x \neq y$, and satisfies $-\int_{\sigma}A \le \cc_A(x,y)$, contradicting $(x,y)\in\I_A$.
        Thus $\ell > 0$, and the unit-speed reparametrization of $\sigma$ is a solution of~\eqref{eq:euler-lagrange-any-speed} with action $\cc_{A}(x,y)$, i.e., a minimizing extremal joining $x$ to $y$.
      \end{proof}

      \begin{lemma}\label{lem:emverification}
        Assume that $M$ is globally hyperbolic.
        Then every precompact set $\cA \subseteq M\times M$ whose closure is contained in $\I_A$ is disjoint from the diagonal and satisfies assumptions~\ref{item:attained} and~\ref{item:compactness}.
      \end{lemma}
      \begin{proof}
        Let $\cA \subseteq M\times M$ be a precompact set whose closure $K$ is contained in $\I_A$.
        Since $M$ is causal and $\cA\subseteq\I_A\subseteq\I$, the set $\cA$ is disjoint from the diagonal.
        In what follows, constants implied by big-O notation may depend also on $K$.
                Let $(x,y) \in K$.
        By~\cite[Th.~1.1]{MinguzziSanchez06}, on each causal homotopy class of causal curves joining $x$ to $y$, the functional $\action^1_A$ attains a minimum at a curve which is either a chronological solution of~\eqref{eq:euler-lagrange-any-speed} or a null geodesic.
        But since $(x,y) \in \I_A$, every null geodesic joining $x$ to $y$ has action greater than $\cc_{A}(x,y)$,
        so we can take a sequence $\gamma_n:[0,\ell_n]\to M$ of
        chronological solutions of~\eqref{eq:euler-lagrange-any-speed}
        joining $x$ to $y$, satisfying $\action_{A}[\gamma_n] \to \cc_{A}(x,y)$ and parametrized by proper time.
        By Lemma~\ref{lem:em-limit}\ref{item:em-limit}, $x$ and $y$ are joined by a minimizing extremal. This proves~\ref{item:attained}.

        \medskip
        It remains to prove~\ref{item:compactness}. Every minimizing extremal $\gamma:[0,\ell]\to M$ joining a pair in $K$ satisfies $|\dot\gamma|_h = O(1/\ell)$ by Lemma~\ref{lem:em-limit}\ref{item:em-speed}. It therefore suffices to bound $\ell$ from below by a positive constant depending only on $K$.

        \medskip
        Let $\gamma_n:[0,\ell_n]\to M$ be a sequence of minimizing extremals with $(\gamma_n(0),\gamma_n(\ell_n))\in K$.
        After passing to a subsequence, we have $(\gamma_n(0),\gamma_n(\ell_n)) \to (x,y) \in K$, and $\limsup_n\action_{A}[\gamma_n] = \limsup_n\cc_{A}(\gamma_n(0),\gamma_n(\ell_n)) \le \cc_{A}(x,y)$ since $\cc_{A}$ is upper semicontinuous by Lemma~\ref{lem:cost-usc}.
        By Lemma~\ref{lem:em-limit}\ref{item:em-limit}, after passing to a further subsequence, $\ell_n$ converges to a positive limit.
        Thus no sequence of minimizing extremals with endpoints in $K$ has proper time
        tending to zero, so the proper time of such extremals is uniformly bounded from below.
      \end{proof}

      \begin{lemma}\label{lem:smallcmr}
        Assume that $M$ is globally hyperbolic.
        Then for every compact set $K \subseteq \I$ there exists $s_0 > 0$ such that $K \subseteq \I_{sA}$ for all $s \in (-s_0,s_0)$.
      \end{lemma}
      \begin{proof}
        Global hyperbolicity implies (see~\cite[Lem.~14.21]{ONeill83})
        that
        \[
          \inf_{K}\tau =:\tau_0 >0
          \qquad
          \text{ and }
          \qquad
          \sup_\gamma
          \Big|\int_\gamma A\Big|
          =: C < \infty
          ,
        \]
        where $\tau$, given by $\tau(x,y):=\sup\int|\dot\gamma|$ with the supremum over causal curves $\gamma$ joining $x$ to $y$, is the time-separation function, and the supremum is taken over causal curves $\gamma$
        joining pairs in $K$.
        Let $(x,y)\in K$ and let $\gamma$ be a maximizing timelike geodesic joining $x$ to
        $y$.
        Then
        \[
          \cc_{sA}(x,y)
          \leq
          \action_{sA}^1[\gamma]
          =
          -
          \tau(x,y)
          -s\int_\gamma A
          \leq
          -\tau_0
          +|s|C
          ,
        \]
        whence, for every null geodesic $\gamma_0$ joining $x$ to $y$,
        \[
          -s\int_{\gamma_0}A
          \geq
          -|s|C
          \geq
          \cc_{sA}(x,y) + \tau_0 - 2|s|C
          .
        \]
        Taking $s_0:=\tau_0/(3C)$, the right-hand side is greater than $\cc_{sA}(x,y)$ for $|s|<s_0$, so that $(x,y)\in\I_{sA}$ by~\eqref{eq:I-A-definition}.
      \end{proof}

    \begin{proof}[Proof of Theorem~\ref{thm:em}.]
      Let $s_0>0$ be given by Lemma~\ref{lem:smallcmr} for the compact set $\supp\mu_0\times\supp\mu_1 \subseteq \I$.
      For $|s|<s_0$ we have $\supp\mu_0\times\supp\mu_1\subseteq\I_{sA}$, so Lemma~\ref{lem:emverification}, applied to the one-form $sA$, enables us to apply Theorem~\ref{thm:diexistence}, which yields a unique monotone displacement interpolation between $\mu_0$ and $\mu_1$.
    \end{proof}

    \bibliographystyle{custom}
    \bibliography{literature}

\end{document}